\documentclass{gtart_h}

\def\ifplaintex{\expandafter\ifx\csname documentclass\endcsname\relax}

\def\gtp{{\mathsurround=0pt\it $\cal G\mskip-2mu$eometry \&\ 
$\cal T\!\!$opology $\cal P\!$ublications}}  

\def\recd{{\small Received:\qua\receiveddate\ifx\reviseddate\relax
\else\qquad Revised:\qua\reviseddate\fi\par}} 

\def\lognumber#1{\def\thelognumber{#1}}
\def\volumenumber#1{\def\thevolumenumber{#1}}
\def\volumeyear#1{\def\thevolumeyear{#1}}
\def\papernumber#1{\def\thepapernumber{#1}}
\def\pagenumbers#1#2{\def\startpage{#1}\def\finishpage{#2}}
\def\published#1{\def\publishdate{#1}}

\def\received#1{\def\receiveddate{#1}}
\def\revised#1{\def\reviseddate{#1}}
\def\accepted#1{\def\accepteddate{#1}}

\let\\\par\let\thelognumber\relax\let\thevolumenumber\relax
\let\thepapernumber\relax\let\thevolumeyear\relax\let\startpage\relax
\let\finishpage\relax\let\publishdate\relax\let\receiveddate\relax
\let\reviseddate\relax\let\accepteddate\relax\let\theasciititle\relax
\let\theasciiauthors\relax
\let\theasciiabstract\relax

\let\theasciiemail\relax

\ifplaintex
\font\logobig=cmssbx10 scaled 3836
\font\logomed=cmssbx10 scaled 2557
\else
\font\logobig=cmssbx10 scaled 4200
\font\logomed=cmssbx10 scaled 2800
\fi

\long\def\makeagttitle{   
\count0=\startpage
\agt\hfill      
\hbox to 45truept{\vbox to 0pt{\vglue -13truept{\logomed A\kern -.37em{\logobig 
T}\kern -.38em G}\vss}\hss}
\break
{\small Volume \thevolumenumber\ (\thevolumeyear)
\startpage--\finishpage\nl
Published: \publishdate}

\vglue .25truein

{\parskip=0pt\leftskip 0pt plus
1fil\def\\{\par\smallskip}{\Large\bf\thetitle}\par\medskip} \vglue
0.05truein

{\parskip=0pt\leftskip 0pt plus 1fil\def\\{\par}{\sc\theauthors}
\par\medskip}%
 
\vglue 0.03truein 

{\small\leftskip 25truept\rightskip 25truept{\bf Abstract}\stdspace\theabstract

{\bf AMS Classification}\stdspace\theprimaryclass
\ifx\thesecondaryclass\relax\else; \thesecondaryclass\fi\par
{\bf Keywords}\stdspace \thekeywords\par}\vglue 7truept

}   

\ifplaintex
\font\phead=cmsl9 scaled 950
\font\pnum=cmbx10 scaled 913
\font\pfoot=cmsl9 scaled 950
\headline{\vbox to 0pt{\vskip -4.5mm\line{\small\phead\ifnum
\count0=\startpage ISSN 1472-2739 (on-line) 1472-2747 (printed)
\hfill {\pnum\folio}\else\ifodd\count0\def\\{ }%
\ifx\theshorttitle\relax\thetitle\else\theshorttitle\fi\hfill{\pnum\folio}
\else\def\\{ and }{\pnum\folio}\hfill\ifx\theshortauthors\relax\theauthors
\else\theshortauthors\fi\fi\fi}\vss}}
\footline{\vbox to 0pt{\vglue 0mm\line{\small\pfoot\ifnum\count0=\startpage
\copyright\ \gtp\hfill\else
\agt, Volume \thevolumenumber\ (\thevolumeyear)\hfill\fi}\vss}}
\else
\font=cmsl9 scaled 1050
\font\lnum=cmbx10 
\font=cmsl9 scaled 1050
\makeatletter
\def\@oddhead{{\small\lhead\ifnum\count0=\startpage ISSN 1472-2739 
(on-line) 1472-2747 (printed)\hfill {\lnum\number\count0}\else\ifodd\count0
\def\\{ }\ifx\theshorttitle\relax \thetitle \else\theshorttitle\fi\hfill
{\lnum\number\count0}\else\def\\{ and }{\lnum\number\count0}
\hfill\ifx\theshortauthors\relax 
\theauthors\else\theshortauthors\fi\fi\fi}}\def\@evenhead{\@oddhead}
\def\@oddfoot{\small\lfoot\ifnum\count0=\startpage\copyright\ \gtp\hfill\else
\agt, Volume \thevolumenumber\ (\thevolumeyear)\hfill\fi}
\def\@evenfoot{\@oddfoot}
\makeatother
\fi
\let\maketitlepage\makeagttitle
\let\makeshorttitle\maketitlepage
\let\maketitle\maketitlepage

\newwrite\gtoutfile
\long\gdef\makeheadfile{  
{\def\\{, }\def\s{ }
\immediate\openout\gtoutfile head.xxx
\immediate\write\gtoutfile{Proxy-for: \ifx\theasciiauthors\relax
\theauthors\else\theasciiauthors\fi\s<\ifx\theasciiemail\relax\theemail\else\theasciiemail\fi>}
\immediate\write\gtoutfile{\noexpand\\}
\immediate\write\gtoutfile{Authors: \ifx\theasciiauthors\relax
\theauthors\else\theasciiauthors\fi}
{\def\\{ }\immediate\write\gtoutfile{Title: \ifx\theasciititle\relax
\thetitle\else\theasciititle\fi}}
\immediate\write\gtoutfile{Subj-class: GT or SG, GR etc}
\immediate\write\gtoutfile{MSC-class: \theprimaryclass\ifx\thesecondaryclass\relax\else, \thesecondaryclass\fi}
\immediate\write\gtoutfile{Journal-ref: Algebr. Geom. Topol. \thevolumenumber\s
(\thevolumeyear) \startpage-\finishpage}
\immediate\write\gtoutfile{Comments: Published by Algebraic and
Geometric Topology at}
\immediate\write\gtoutfile{\s\s\s  http://www.maths.warwick.ac.uk/agt/AGTVol\thevolumenumber/agt-\thevolumenumber-\thepapernumber.abs.html}
\immediate\write\gtoutfile{\noexpand\\}
\immediate\write\gtoutfile{}
\ifx\theasciiabstract\relax
\immediate\write\gtoutfile{\theabstract}\else
\immediate\write\gtoutfile{\theasciiabstract}\fi
\immediate\write\gtoutfile{}
\immediate\write\gtoutfile{\noexpand\\}
\immediate\write\gtoutfile{}
\immediate\closeout\gtoutfile}}  

\def\maketitlepage{\makeagttitle\makeheadfile}
\let\makeshorttitle\maketitlepage
\let\maketitle\maketitlepage

\lognumber{25}
\volumenumber{5}
\volumeyear{2005}
\papernumber{25}
\pagenumbers{577}{613}
\received{13 February 2005} 
\revised{28 April 2005}
\accepted{31 May 2005}
\published{30 June 2005}

\usepackage{amsmath, amssymb, amsthm}
\usepackage{labelfig, graphicx}
\let\Bbb\mathbb
\theoremstyle{plain}
\newtheorem{thm}{Theorem}[section]
\newtheorem{lem}[thm]{Lemma}
\newtheorem{cor}[thm]{Corollary}
\newtheorem{prop}[thm]{Proposition}
\newcommand{\oddspin}[1]{\mathcal{SP}_{#1} [q_1]}
\newcommand{\spin}[2]{\mathcal{SP}_{#1}[#2]}
\newcommand{\Diff}{\mathrm{Diff}}
\newcommand{\inter}{\mathrm{int}}
\newcommand{\fix}{\mathrm{fix}}
\newcommand{\Ortho}{\mathrm{O}}
\newcommand{\Aut}{\mathrm{Aut}}
\newcommand{\Symp}{\mathrm{Sp}}
\newcommand{\inv}[1]{\overline{#1}\;}

\newcommand{\Ggeq}{\underset{G_g}{\sim}}

\newcommand{\Cbullet}{\bullet\bullet\bullet}
\newcommand{\repre}[1]{\overline{\overline{#1}}}
\newcommand{\Jump}{\text{Lemma \ref{lem:Gg-action}(1)}}
\newcommand{\Shift}{\text{Lemma \ref{lem:Gg-action}(2)}}
\newcommand{\complex}{{\mathbb{C}}}
\newcommand{\ComplexPlane}{{\mathbb{CP}}^2}
\def\co{\colon\thinspace}
\newfont{\ttnonit}{cmtt10}
\theoremstyle{definition}

\theoremstyle{remark}
\newtheorem{rem}[thm]{Remark}
\newtheorem{prf}{Proof}

\begin{document}
\title{Surfaces in the complex projective plane\\and their mapping class groups}
\author{Susumu Hirose}
\address{Department of Mathematics, 
Faculty of Science and Engineering\\Saga University, 
Saga, 840-8502 Japan.}
\email{hirose@ms.saga-u.ac.jp}
\begin{abstract} 
An orientation preserving diffeomorphism over a surface embedded in 
a 4-manifold is called extendable, 
if this diffeomorphism is a restriction of an orientation preserving 
diffeomorphism on this 4-manifold. 
In this paper, we investigate conditions for extendability of diffeomorphisms over 
surfaces in the complex projective plane. 
\end{abstract}
\maketitle
\cl{\small\it Dedicated to Professor Yukio Matsumoto for his 60th birthday} 
\section{Introduction}\label{sec:intro}
%
There are deformations of embedded surfaces in 4-manifolds which induce 
isotopically non-trivial diffeomorphisms on surfaces. 
We introduce two typical examples. 

\begin{figure}[ht!]\small
\begin{center}
\SetLabels
(0.15*0.5) (1) \\
(0.49*0.5) (2) \\
(0.85*0.5) (3) \\
(0.15*-0.05) (4) \\
(0.49*-0.05) (5) \\
(0.85*-0.05) (6) \\
(0.92*0.8) $A$ \\
(0.22*0.22) $A$ \\
\endSetLabels
\strut\AffixLabels{\includegraphics[height=5cm]{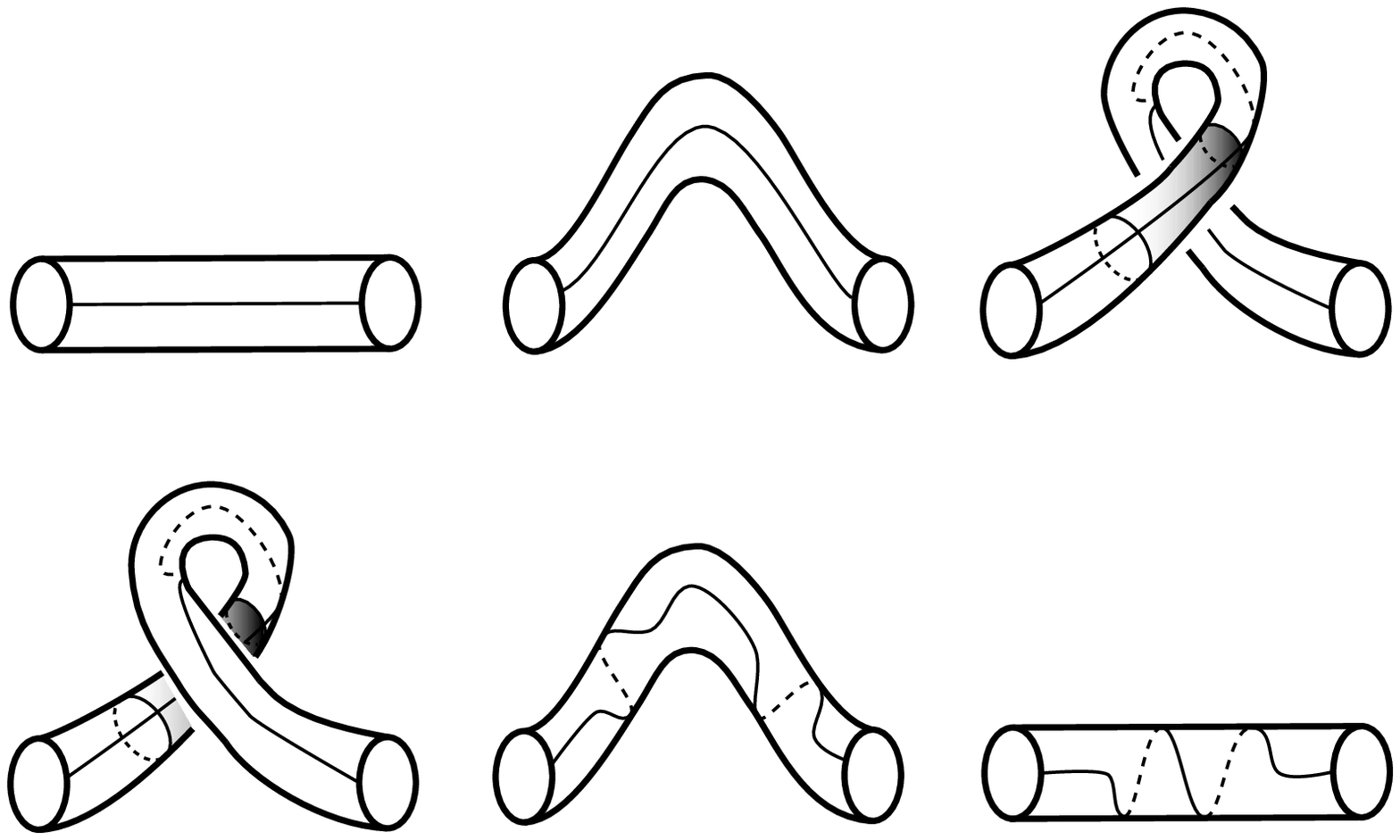}}
\vspace{2mm}\nocolon
\caption{}
\label{fig:bandmove}
\end{center}
\end{figure}
For the first example, we consider a deformation of an annulus embedded in 
$S^3 \times [-1, 1]$ so that, under this deformation, the boundary of this annulus 
is fixed. 
Let $S^1 \times [0,1]$ be an annulus embedded in $S^3 \times \{ 0 \}$ 
$ \subset S^3 \times [-1, 1]$, and 
$t \co S^3 \times [-1,1] \to [-1,1]$ a projection to the second factor. 
We deform $S^1 \times [0,1]$ as in Figure \ref{fig:bandmove}. 
First, we isotope $S^1 \times [0,1]$ in $S^3$ from (1) to (3). 
Next, we isotope $S^1 \times [0,1]$ so that outside of the annulus $A$ of (3) 
$t=0$, and inside $t > 0$. 
Then we isotope $S^1 \times [0,1]$ inside $A$ so that, when we push $A$ down to 
$S^3 \times \{ 0 \}$, $S^1 \times [0,1]$ is as in (4). 
Finally, we isotope $S^1 \times [0,1]$ in $S^3$ from (4) to (6). 
The composition of these deformations induce a square of Dehn twist about 
the core circle $S^1 \times \{ \frac{1}{2} \}$ of $S^1 \times [0,1]$. 

For the second example, we consider a deformation of a non-singular plane curve 
of degree $3$. A torus is defined as a quotient of the complex plane by a lattice 
${\mathbb Z} + {\mathbb Z} \sqrt{-1}$. 
We embed this torus into the complex projective plane by using the Weierstrass 
$\wp$ function associated to this lattice, then the image of this embedding is 
a non-singular plane curve of degree $3$. 
We deform this lattice, ${\mathbb Z} + {\mathbb Z} (\sqrt{-1} + t)$, 
where $0 \leq t \leq 1$ is a parameter of this deformation. 
Then the embedding is deformed isotopically and, finally (when $t=1$), 
brought back to the original position. 
This deformation induces a Dehn twist on the non-singular plane curve of degree $3$. 

In this paper, we investigate a topological meaning of the above phenomena. 

\medskip

We settle a general formulation. 
Let $M$ be a simply connected compact oriented smooth 4-manifold 
(possibly with boundary) 
and $F$ be a compact oriented smooth 2-manifold (possibly with boundary) 
embedded in $M$. 
We call the pair $(M, F)$ a {\em knotted surface\/}. 
In particular, if $F$ is characteristic, that is to say, 
$F \cdot X \equiv X \cdot X \mod 2$ 
for any $X \in H_2 (M, \mathbb{Z})$, then we call 
this pair $(M, F)$ a {\em knotted characteristic surface\/}. 
An orientation preserving diffeomorphism $\psi$ over $F$ is 
{\em extendable\/} if there is an orientation preserving diffeomorphism 
$\Psi$ over $M$ such that $\Psi |_{F} = \psi$. 
In general, for an oriented manifold $A$ and its submanifold $B$, 
we denote 
$$
\Diff_+ (A, \fix \, B) = \left\{ \psi \; \left| \;
\begin{aligned}
&\text{ an orientation preserving 
diffeomorphism over }A \\
&\text{ such that } \psi|_{B} = id_B 
\end{aligned}
\right.
\right\}. 
$$ 
If $B=\phi$, we denote this group by $\Diff_+(M)$. 
The group $\pi_0 (\Diff_+ (F, \fix \, \partial F))$ is called the 
{\em mapping class group\/} of $F$ and denoted by ${\mathcal{M}}_{F}$. 
If $F$ is a closed oriented surface of genus $g$, this group 
is denoted by ${\mathcal{M}}_g$. 
We define 
$$
{\mathcal{E}}(M, F) = \{\psi \in {\mathcal{M}}_{F} \; | \; 
\psi \text{ is extendable } \}. 
$$
This is a subgroup of ${\mathcal{M}}_{F}$ and 
is a central object of this paper. 

In the case where $M = S^4$, there are several works on this group. 
Let $(S^4, \Sigma_g)$ be the genus $g$ trivial knotted surface in 
$S^4$. 
When $g=1$, Montesinos \cite{Montesinos} investigated 
${\mathcal{E}}(S^4, \Sigma_1)$, and when $g \geq 2$, 
the author \cite{Hirose2} investigated ${\mathcal{E}}(S^4, \Sigma_g)$. 
Let $(S^3,k)$ be a knot in $S^3$ and $(S^4, S(k))$ (resp. 
$(S^4, \tilde{S}(k))$) the spun (resp. the twisted spun) of 
$(S^3, k)$. 
When $(S^3, k)$ is a torus knot, Iwase \cite{Iwase} investigated 
${\mathcal{E}}(S^4, S(k))$ and ${\mathcal{E}}(S^4, \tilde{S}(k))$, 
and when $(S^3, k)$ is an arbitrary knot, 
the author \cite{Hirose} investigated these groups. 

In this paper, we investigate the case where $M$ is 
the complex projective plane $\ComplexPlane$. 
In \S \ref{sec:standard}, we treat the case where $(\ComplexPlane, \Sigma_g)$ is 
a standard embedding of $\Sigma_g$. In \S \ref{sec:curve}, we treat the case 
where $(\ComplexPlane, F)$ is a non-singular plane curve. 
From \S \ref{sec:connected-sum} to the end of this paper, we treat the case 
where $(\ComplexPlane, F)$ is a connected sum of a non-singular plane curve of 
degree $3$ and a trivial embedding. 

\section{Preliminary: A Hopf band on the boundary of the 4-ball}
\label{sec:preliminary}
%
%
A link $L$ in $S^3$ is called a {\em fibered link\/} if there is a map 
$\phi \co S^3 \setminus L \to S^1$ which is a fiber bundle projection. 
For each $t \in S^1$, $\phi^{-1} (t) = F$, which does not depend on $t$, 
is called the {\em fiber\/} of $L$. 
Since $\phi$ is a bundle projection, $S^3 \setminus L$ is diffeomorphic to 
the quotient of $F \times [0,1]$ by an equivalence 
$(x,0) \sim (h(x), 1)$ where $h$ is a diffeomorphism over $F$ and called 
the {\em monodromy\/} of $L$. 

\begin{figure}[ht!]\small
\begin{center}
\SetLabels
\L(0.04*-0.18) Positive Hopf band \\
\L(0.56*-0.18) Negative Hopf band \\
(0.36*0.8) $c$ \\
(0.885*0.8) $c$ \\
\endSetLabels
\strut\AffixLabels{\includegraphics[height=2cm]{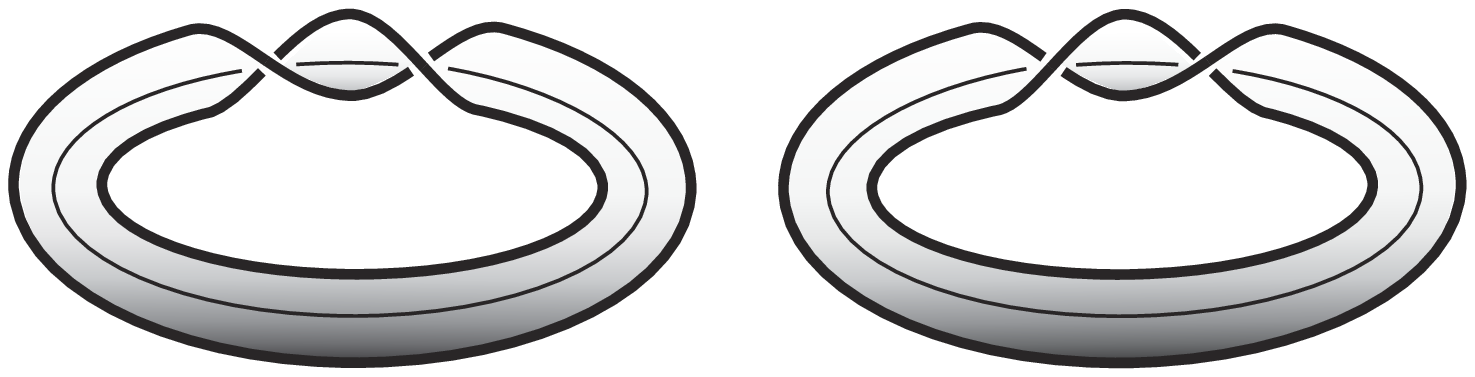}}
\vspace{2mm}
\nocolon
\caption{}
\label{fig:Hopf-band}
\end{center}
\end{figure}

A {\em Hopf band\/} is an annulus embedded in $S^3$ as in Figure 
\ref{fig:Hopf-band}. In this picture, there are two types of Hopf bands. 
In this note, we treat both types of Hopf bands. 
The boundary of a Hopf band is called a {\em Hopf link\/}. 
The Hopf link is a fibered link whose fiber is the Hopf band and 
whose monodromy is a Dehn twist about the core circle of the Hopf band. 
Let $B$ be a Hopf band in $S^3$ which is a boundary of a 4-ball $D^4$. 
We push the interior of $B$ into the interior of $D^4$ and 
let $B'$ be the annulus 
obtained by this deformation and let $c$ be the core circle of 
$B'$. 
\begin{prop}\label{prop:Hopf-twist}
The Dehn twist $T_c$ about $c$ is extendable, 
i.e. there is an element $T$ $\in$ 
$\Diff_+ (D^4, \fix \, \partial D^4)$ 
such that $T|_{B'} = T_c$. 
\end{prop}
\begin{prf}
Since $\partial B$ is a fibered link, whose fiber is $B$ and 
whose monodromy is $T_c$, there is an orientation preserving diffeomorphism 
$\psi$ of $S^3$ such that $\psi|_{B} = T_c$, 
and there is an isotopy $\psi_t$ ($t \in [0,1]$) with $\psi_0 = id_{S^3}$ 
and $\psi_1 = \psi$, which is defined by shifting fibers. 
Let $N(\partial D^4)$ be the regular neighborhood of $\partial D^4$ in $D^4$. 
We parametrize $N(\partial D^4) = S^3 \times [0,2]$ so that 
$S^3 \times \{0 \} = \partial D^4$ and $B' = \partial B \times [0,1] \cup 
B\times \{1\}$. 
Let $T$ be a diffeomorphism defined as follows 
\begin{align*}
&T|_{N(\partial D^4)}(x,t) = 
\begin{cases}
(\psi_t(x),t) & 0 \leq t \leq 1 \\
(\psi_{2-t}(x), t) & 1 \leq t \leq 2 
\end{cases} \\
&T|_{D^4 \setminus N(\partial D^4)} = id. 
\end{align*}
This is the diffeomorphism which we need. 
\qed
\end{prf}
\begin{rem} Let $(S^4, \Sigma_g)$ be the genus $g$ surface 
standardly embedded in $S^4$. 
In \cite{Hirose2}, the author showed that 
$T_{c_4} T_{c_3} T_{c_4}^{-1}$ $\in$ ${\mathcal{E}}(S^4, \Sigma_g)$ 
by using Montesinos' result \cite[Theorem 5.3]{Montesinos}
($c_3$ and $c_4$ are as in Figure \ref{fig:circle}) . 
We show this fact by using Proposition \ref{prop:Hopf-twist}. 
The 4-sphere $S^4$ is constructed from two 4-balls $D^4_+$, $D^4_-$ 
with attaching along the boundary $S^3 = \partial D^4_+ = \partial D^4_-$. 
We parametrize the regular neighborhood $N(\partial D^4_+) = S^3 \times [0,2]$ 
in $D^4_+$ so that $\partial D^4_+ = S^3 \times \{ 0 \}$. 
The regular neighborhood $N$ of $T_{c_4} (c_3)$ in $\Sigma_g$ is 
a Hopf band in $S^3$ $\subset S^4$. 
We push the interior of $N$ into the interior of $D^4_+$, then
we get an annulus $N'$ properly embedded in $D^4_+$. 
We may assume, by the above parametrization of $N(\partial D^4_+)$, 
$N' \cap S^3 \times \{ t \} = \partial N \times \{ t \}$ for $0 \leq t < 2$ 
and $N' \cap S^3 \times \{ 2 \} = N \times \{ 2 \}$. 
We denote $D^4_+ \setminus S^3 \times [0,1)$ by $D'$. 
By applying Proposition \ref{prop:Hopf-twist} to $(D', N' \cap D')$, 
we show that there is an element $T$ $\in$ $\Diff_+ (D', 
\fix \, \partial D')$ such that $T|_{N' \cap D'} = T_{c_4} T_{c_3} T_{c_4}^{-1}$. 
Therefore, we see $T_{c_4} T_{c_3} T_{c_4}^{-1}$ $\in$ 
${\mathcal{E}}(S^4, \Sigma_g)$. 
\end{rem}
\section{Surfaces standardly embedded in the complex projective plane}
\label{sec:standard}
%
For the free action of ${\mathbb{C}}^* = {\mathbb{C}}\setminus \{ 0 \}$ on 
${\mathbb{C}}^3 \setminus \{ (0,0,0) \}$ defined by 
$\lambda (z_0, z_1 , z_2) = (\lambda z_0, \lambda z_1, \lambda z_2)$, 
we take the quotient 
$\ComplexPlane = (\complex^3 \setminus \{ (0,0,0) \})/ {\mathbb{C}}^*$. 
This space $\ComplexPlane$ is a closed oriented 4-manifold and 
called the {\em complex projective plane\/}. 
This 4-manifold $\ComplexPlane$ is constructed from $D^4$ by 
attaching a 2-handle $h^2$ along the frame $1$ trivial knot $K_0$ in 
$\partial D^4$, and attaching a 4-handle $h^4$. 
A {\it 3-dimensional handlebody\/} $H_g$ is an oriented 3-manifold 
which is constructed from a 3-ball with attaching $g$ 1-handles. 
Any image of embeddings of $H_g$ into $\ComplexPlane$ are isotopic each other. 
Therefore, $(\ComplexPlane, \partial H_g)$ is unique. 
A surface {\it standardly embedded} in $\ComplexPlane$ is 
$(\ComplexPlane, \partial H_g)$. We obtain: 
\begin{thm}\label{thm:standard}
For any $g$, ${\mathcal{E}}(\ComplexPlane, \partial H_g)= {\mathcal{M}}_g$. 
\end{thm}
\begin{prf}
Let $D^4$ be the 4-ball used to construct $\ComplexPlane$ and 
$N(\partial D^4)$ be the regular neighborhood of $\partial D^4$ in $D^4$. 
We parametrize $N(\partial D^4) = S^3 \times [0,-1]$, 
so that $S^3 \times \{ 0 \} = \partial D^4$ and, for $-1 \leq t < 0$, 
$S^3 \times \{ t \}$ is in the interior of $D^4$. 
Since the image of embedding of $H_g$ in $\ComplexPlane$ is unique up to 
isotopy, we assume that $H_g$ $\subset$ $S^3 \times \{ -1 \}$ 
and that each simple closed curve $c$ on $H_g$ which corresponds to Lickorish 
generator of mapping class group ${\mathcal{M}}_g$ is a trivial knot in 
$S^3 \times \{ -1 \}$. 
The regular neighborhood $N(c)$ of $c$ on $\partial H_g$ 
is an annulus trivially embedded in $S^3 \times \{ -1 \}$. 
At first, we deform $H_g$ in $S^3 \times \{ -1 \}$ 
so that, if we forget the second factor $[0,-1]$, 
$c \cup K_0$ becomes a Hopf link in $S^3$. 
We push $N(c)$ into $\partial (D^4 \cup h^2)$, then 
$N(c)$ becomes a Hopf band in $\partial h^4$. 
By applying Proposition \ref{prop:Hopf-twist}, 
we see that $T_{c}$ is extendable in $h^4$, and so in $\ComplexPlane$. 
\qed
\end{prf}

\section{Non-singular plane curves}\label{sec:curve}
%
%
We review here the topological description of non-singular plane curves 
by Akbulut and Kirby \cite{A-K} (see also \cite[6.2.7]{G-S}). 

\begin{figure}[ht!]\small
\begin{center}
\SetLabels
(0.35*0.1) $F_{3,3}$ \\
(0.21*0.5) $c_1$ \\ (0.26*0.7) $c_2$ \\
(0.97*0.1) $F_{4,4}$ \\
(0.52*0.6) $c_1$ \\ (0.67*0.65) $c_2$ \\ (0.69*0.34) $c_3$ \\ 
(0.75*0.45) $c_4$ \\ (0.72*0.58) $c_5$ \\ (0.7*0.15) $d$ \\
(0.88*0.65) $c_6$ \\ (0.8*0.08) $c_7$ \\
\endSetLabels
\strut\AffixLabels{\includegraphics[height=6cm]{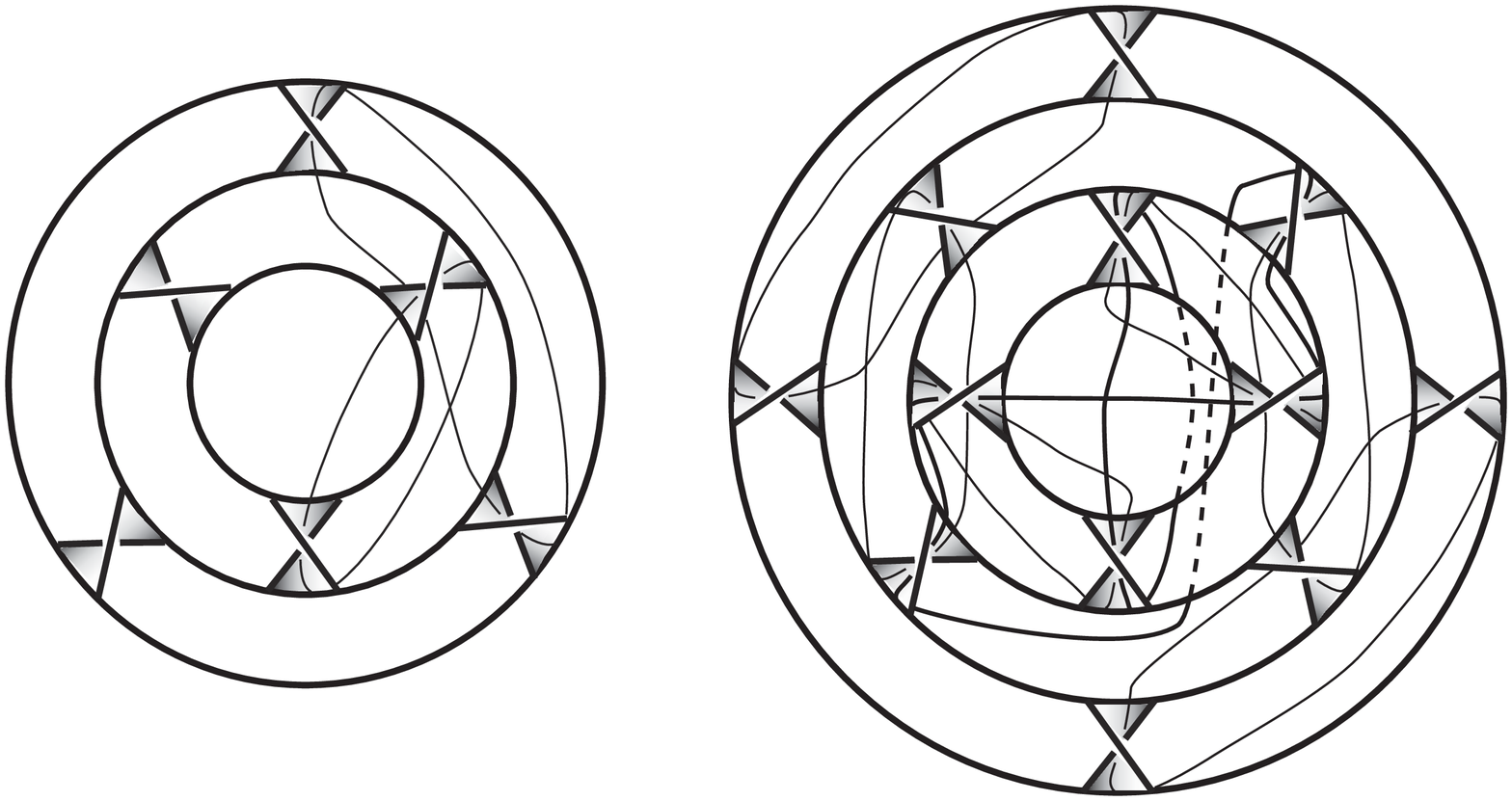}}
\nocolon
\caption{}
\label{fig:Seifert-surface-12}
\end{center}
\end{figure}

\begin{figure}[ht!]\small
\begin{center}
\SetLabels
\R(0.01*0.8) $T_{d,d}$ \\ (0.37*0.42) $K_0$ \\
\L(0.38*0.3) attach a $2$-handle \\ \L(0.38*0.22) along $K_0$ \\
\endSetLabels
\strut\AffixLabels{\includegraphics[height=5cm]{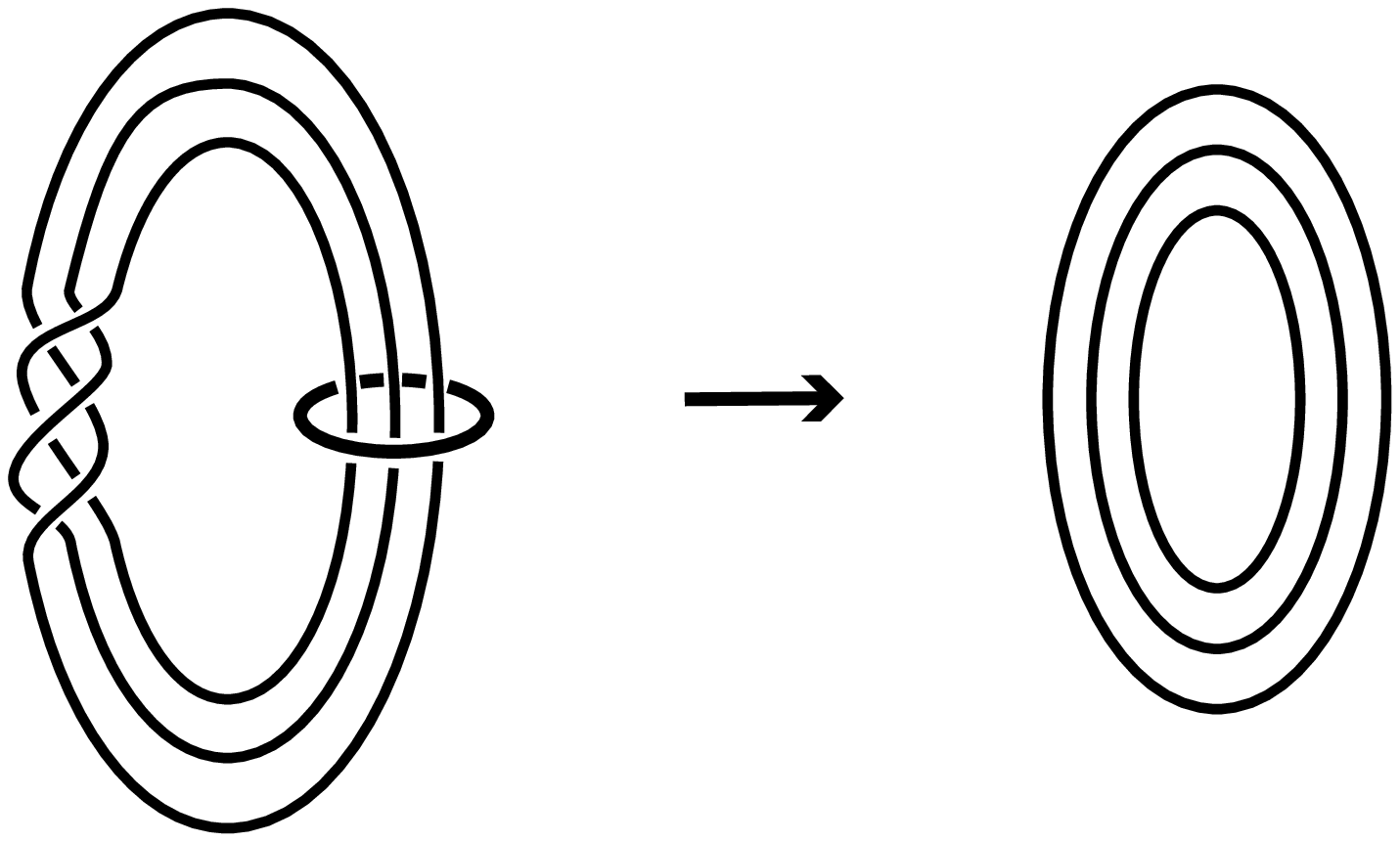}}
\nocolon
\caption{}
\label{fig:C3}
\end{center}
\end{figure}

An $(m,n)$-{\em torus link} $T_{m,n}$ is an oriented link in 
$S^3 = \partial D^4$ consisting of $\gcd(m,n)$ oriented circles in 
the boundary of the tubular neighborhood $U$ of the trivial knot, 
representing $m \mu + n \lambda$ in $H_1(\partial U ; {\mathbb{Z}})$, 
where $\mu = [$the meridian of the trivial knot$]$ and 
 $\lambda = [$the longitude of the trivial knot$]$. 
There is a canonical Seifert surface $F_{m,n}$ for $T_{m,n}$, 
consisting of $n$-disks connected by $m(n-1)$ twisted bands as 
in Figure \ref{fig:Seifert-surface-12}. 
As $K_0$, we take a trivial knot given by pushing $T_{1,0}$ into 
the complement of $U$ ( see the left hand side of 
Figure \ref{fig:C3}). 
From here, we consider only the case where $m=n=d$. 
As shown in the right hand side of Figure \ref{fig:C3}, 
$T_{d,d}$ becomes $d$ components trivial link in 
$\partial (D^4 \cup h^2)$. 
Let $D_d$ be disjoint 2-disks in $\partial (D^4 \cup h^2)$ 
which bound this trivial link. 

Let $K_d$ be a non-singular plane curve of degree $d$, then $K_d$ is a 
genus $\frac{(d-1)(d-2)}{2}$ closed oriented surface embedded in 
$\ComplexPlane$. 
We remark that $K_d$ is unique up to isotopy, 
$K_d =\{ [X:Y:Z] \in \ComplexPlane | X^d + Y^d + Z^d = 0 \}$ 
and $[K_d] = d [{\mathbb{CP}}^1]$ 
$\in H_2(\ComplexPlane; {\mathbb{Z}})$. 
Akbulut and Kirby showed: 
\begin{prop}\label{prop:Akbulut-Kirby}
$K_d = F_{d,d} \cup D_d. $ 
\end{prop}
Thus we obtain: 
\begin{thm}\label{thm:non-sing-low}
When $d=3,4$, ${\mathcal{E}}(\ComplexPlane, K_d)= {\mathcal{M}}_{g_d}$, 
where $g_d = \frac{(d-1)(d-2)}{2}$. 
\end{thm}
\begin{prf}
When $d=3$, $K_3$ is homeomorphic to a 2-dimensional torus $T^2$. 
In $F_{3,3}$ (see Figure \ref{fig:Seifert-surface-12}), 
each regular neighborhood of $c_1$ and $c_2$ is a Hopf band. 
Therefore, by Proposition \ref{prop:Hopf-twist}, $T_{c_1}$ and $T_{c_2}$ 
are elements of ${\mathcal{E}}(\ComplexPlane, K_3)$. 
On the other hand, $T_{c_1}$ and $T_{c_2}$ generate ${\mathcal{M}}_1$. 
Hence, ${\mathcal{E}}(\ComplexPlane, K_3) = {\mathcal{M}}_1$. 

\begin{figure}[ht!]\small
\begin{center}
\SetLabels
(0.07*0.35) $c_1$ \\ (0.35*0.35) $c_3$ \\ (0.66*0.35) $c_5$ \\ 
(0.94*0.35) $c_7$ \\ (0.2*0.2) $c_2$ \\ (0.5*0.18) $c_4$ \\
(0.81*0.19) $c_6$ \\ (0.55*0.82) $d$ \\
\endSetLabels
\strut\AffixLabels{\includegraphics[height=2.5cm]{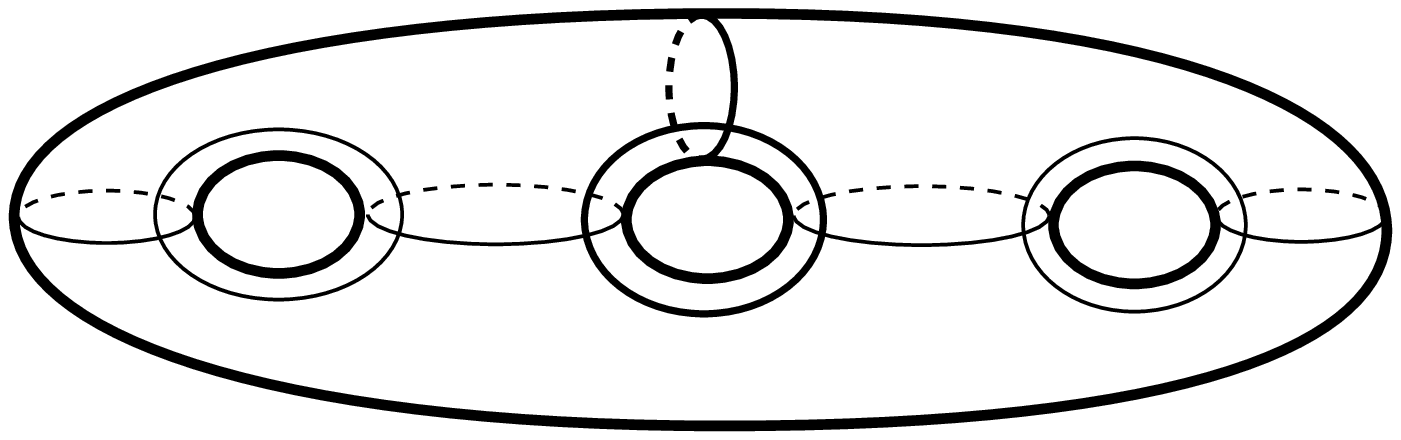}}
\nocolon
\caption{}
\label{fig:Lickorish-generator-3}
\end{center}
\end{figure}

When $d=4$, we do the same as the above case. We remark that the Dehn twists about 
the simple closed curves in Figure \ref{fig:Lickorish-generator-3} 
corresponding to the simple closed curves 
in $F_{4,4}$ (see Figure \ref{fig:Seifert-surface-12}) with the same symbols 
generate the mapping class group of genus $3$ surface 
\cite{Lickorish}. 
\qed
\end{prf}
When $d \geq 5$, ${\mathcal{E}}(\ComplexPlane, K_d)$ is unknown. 
It is, however, not the case that 
${\mathcal{E}}(\ComplexPlane, K_d)$ $=$ ${\mathcal{M}}_{g_d}$, 
because, when $d$ is odd, $K_d$ is a characteristic surface, 
so the Rokhlin quadratic form on $H_1(K_d; {\mathbb{Z}}_2)$ 
is well-defined (we review the definition of the Rokhlin quadratic form 
in the next section). 
By the definition of the Rokhlin quadratic form, if a diffeomorphism 
on $K_d$ is extendable to $\ComplexPlane$, this diffeomorphism should 
preserve this form. Hence: 
\begin{thm}\label{thm:non-sing-high-odd}
When $d$ is an odd integer greater than or equal to $5$, 
${\mathcal{E}}(\ComplexPlane, K_d)$ is a proper subgroup of 
${\mathcal{M}}_{g_d}$, 
where $g_d = \frac{(d-1)(d-2)}{2}$. 
\end{thm}

\section{Connected sum of the non-singular plane curve of degree $3$ and 
trivial knotted surface}
\label{sec:connected-sum}
%

We define knotted surfaces investigated from here to the end of this paper. 
The images of any embeddings of a $3$-dimensional handlebody $H_g$ 
into $S^4$ are isotopic each other. 
We call this $\Sigma_g$-knot $(S^4, \partial H_g)$ 
{\it a trivial $\Sigma_g$-knot\/}, and this is denoted by $(S^4, \Sigma_g)$. 
Let $(\ComplexPlane, K_3)$ be a nonsingular cubic plane curve. 
We define connected sum of $(\ComplexPlane, K_3)$ and $(S^4, \Sigma_{g-1})$ 
following the construction by Boyle \cite{Boyle} as follows. 
We choose points $p$ and $q$ on $K_3$ and $\Sigma_{g-1}$ respectively, and find 
small 4-balls $B_1$ and $B_2$ centered at $p$ and $q$ such that 
the pairs $(B_1, B_1 \cap K_3)$ and $(B_2, B_2 \cap \Sigma_{g-1})$ are equivalent 
to the standard pair $(B^4, B^2)$. Now we glue the pairs 
$(S^4 \setminus \inter(B_1), K_3 \setminus \inter(B_1))$ and 
$(\ComplexPlane \setminus \inter(B_2), \Sigma_{g-1} \setminus \inter(B_2))$ 
together by an orientation-reversing diffeomorphism 
$f:\partial B_1 \to \partial B_2$ such that 
$f(\partial B_1 \cap K_3) = \partial B_2 \cap \Sigma_{g-1}$. 
Since the connected sum of $\ComplexPlane$ and $S^4$ is diffeomorphic to 
$\ComplexPlane$, we get a surface in $\ComplexPlane$ and denote this 
characteristic knotted surface by $(\ComplexPlane, K_3 \# \Sigma_{g-1})$. 
From here to the end of this paper, we investigate on the group 
${\mathcal{E}}(\ComplexPlane, K_3 \# \Sigma_{g-1})$. 

For a knotted characteristic surface $(M, F)$, where $M$ is a simply connected 
smooth closed oriented 4-manifold, 
we define a quadratic form ({\it the Rokhlin quadratic form\/}) 
$q_F: H_1(F; {\Bbb Z}_2) \to {\Bbb Z}_2$: 
Let $P$ be a compact surface embedded in $M$, 
with its boundary contained in $F$, normal to $F$ along its 
boundary, and its interior is transverse to $F$. 
Let $P'$ be a surface transverse to $P$ obtained by sliding $P$ parallel to 
itself over $F$. 
Define $q_F([\partial P]) = \#( \text{int} P \cap (P' \cup F)) \text{
mod }  2$. 
This is a well-defined quadratic form with respect to 
the ${\Bbb Z}_2$-homology intersection form $(,)_2$ on $F$, 
i.e. for each pair of elements $x$, $y$ of 
$H_1(F; {\Bbb Z}_2)$, $q_F(x+y) = q_F(x) + q_F(y) + (x,y)_2$. 
By the definition of the Rokhlin quadratic from $q_F$, 
if $\psi \in \Diff_+ (F)$ is extendable, then $\psi$ preserves $q_F$, 
that is to say, $q_F (\psi_*(x)) = q_F(x)$ for any $x$ $\in$ 
$H_1 (F; {\Bbb Z}_2)$. 
We will show,  

\begin{thm}\label{thm:main}
%
For any $g \geq 2$, 
$$
{\mathcal E}(\ComplexPlane, K_3 \# \Sigma_{g-1}) = 
\left\{ \psi \in {\mathcal{M}}_g \; \left| \;
\begin{aligned}
&
q_{K_3 \# \Sigma_{g-1}} (\psi_*(x)) = q_{K_3 \# \Sigma_{g-1}} (x) \\
&\text{ for any } x \in H_1 (K_3 \# \Sigma_{g-1}; {\Bbb Z}_2)
\end{aligned}
\right.
\right\} .   
$$
\end{thm}

In \S \ref{sec:generators}, we investigate on a system of generators for 
the right hand side group in the equation of Theorem \ref{thm:main}. 
In \S \ref{sec:extension}, we show that each element of this system of generators 
is extendable. 

\section{A finite set of generators for the odd spin mapping class group}
\label{sec:generators}
%
%
\begin{figure}[ht!]\small
\begin{center}
\SetLabels
(0.05*0.4) $\gamma_1$ \\ (0.145*0.48) $\gamma_2$ \\ (0.22*0.4) $\gamma_3$ \\
(0.3*0.48) $\gamma_4$ \\ (0.33*0.16) $\beta'_4$ \\ (0.33*0.82) $\beta_4$ \\
(0.38*0.4) $\gamma_5$ \\
(0.45*0.48) $\gamma_6$ \\ (0.48*0.16) $\beta'_6$ \\ (0.48*0.82) $\beta_6$ \\
(0.7*0.48) $\gamma_{2g-2}$ \\ (0.75*0.16) $\beta'_{2g-2}$ \\ 
(0.75*0.82) $\beta_{2g-2}$ \\
(0.78*0.4) $\gamma_{2g-1}$ \\ (0.86*0.48) $\gamma_{2g}$ \\ 
(0.94*0.4) $\gamma_{2g+1}$ \\ (0.03*0.05) $\gamma_{0}$ \\
\endSetLabels
\strut\AffixLabels{\includegraphics[height=4cm]{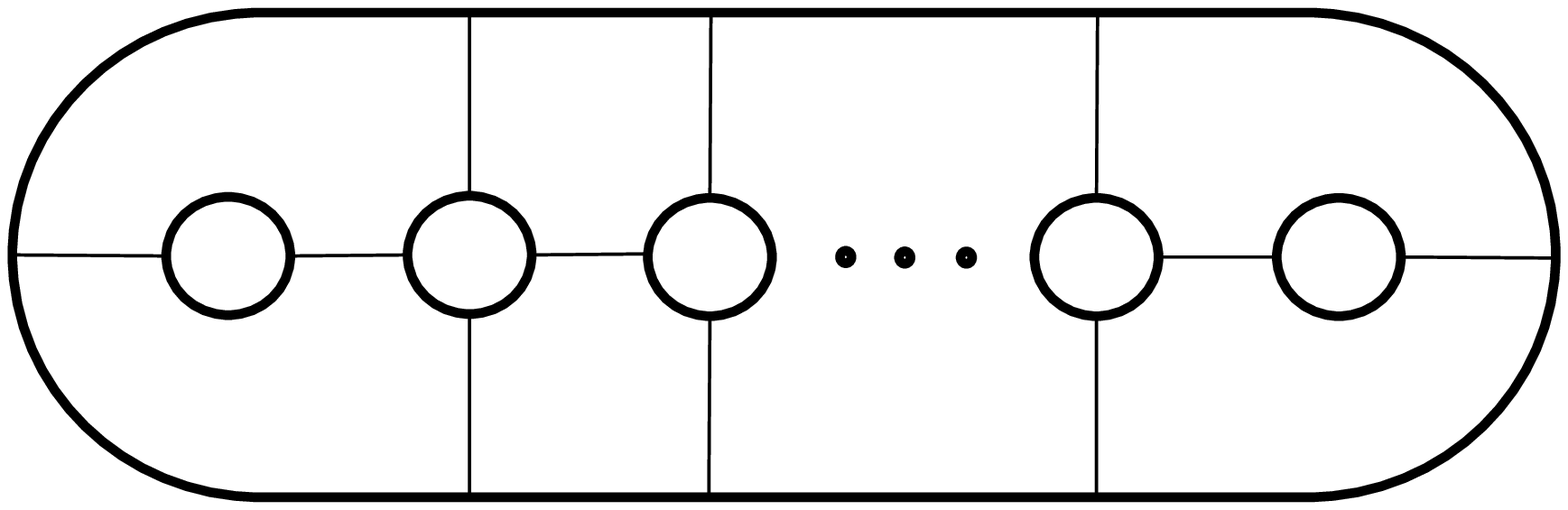}}
\nocolon
\caption{}
\label{fig:planar}
\end{center}
\end{figure}
\begin{figure}[ht!]\small
\centering
\SetLabels
(0.04*0.32) $c_1$ \\ (0.115*0.21) $c_2$ \\ (0.19*0.30) $c_3$ \\ 
(0.27*0.21) $c_4$ \\ (0.35*0.8) $c_{\beta}$ \\ (0.72*0.21) $c_{2g-2}$ \\
(0.79*0.30) $c_{2g-1}$ \\ (0.87*0.21) $c_{2g}$ \\ \L(1.01*0.49) $c_{2g+1}$ \\
\endSetLabels
\strut\AffixLabels{\includegraphics[height=3.5cm]{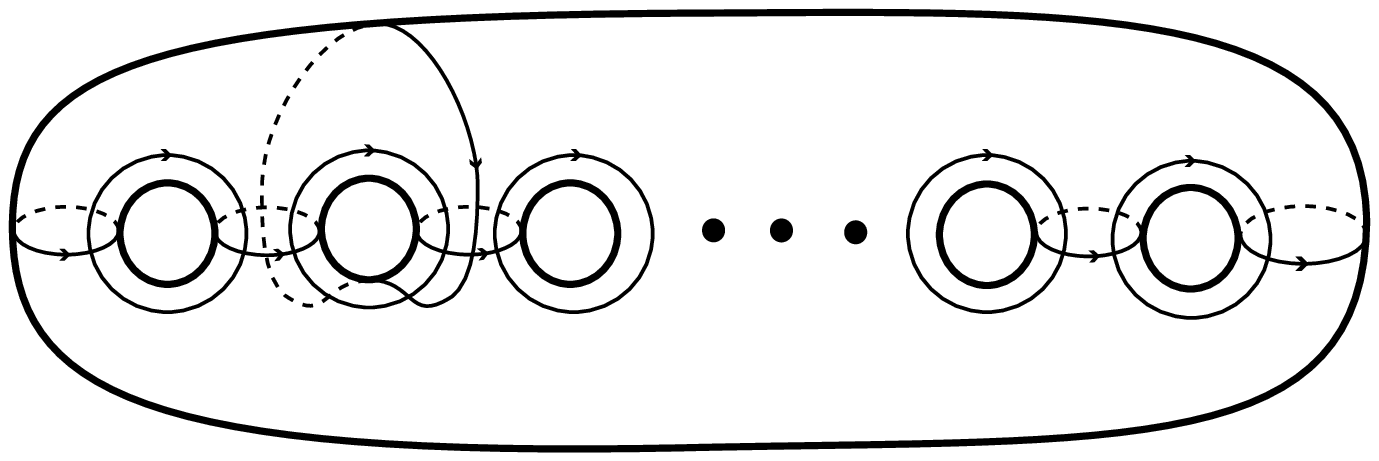}}
\nocolon
\caption{}
\label{fig:circle}
\end{figure}
\begin{figure}[ht!]\small
\centering
\SetLabels
(0.3*0.76) $b_4$ \\ (0.3*0.19) $b_4'$ \\ (0.44*0.76) $b_6$ \\ 
(0.44*0.19) $b_6'$ \\ (0.73*0.76) $b_n$ \\ (0.73*0.19) $b_n'$ \\
(0.75*0.35) $c_n$ \\
\endSetLabels
\strut\AffixLabels{\includegraphics[height=3cm]{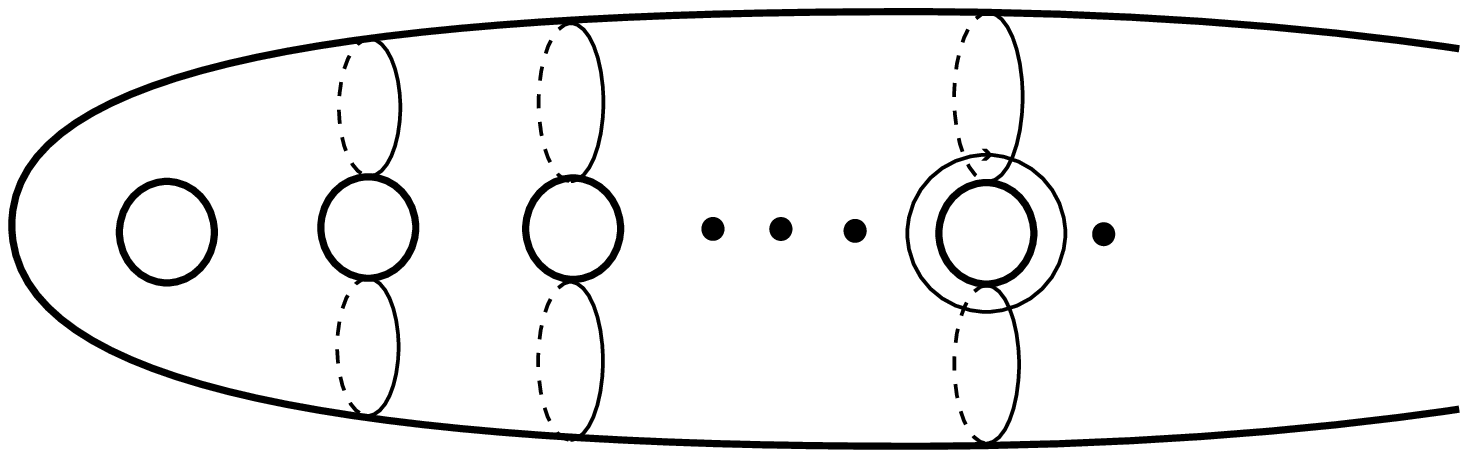}}
\nocolon
\caption{}
\label{fig:b-circle}
\end{figure}

We settle some notations. 
Let $P_g$ be a planar surface constructed from a 2-disk by removing 
$g$ copies of disjoint 2-disks. 
As indicated in Figure \ref{fig:planar}, we denote the boundary components of 
$P_g$ by $\gamma_0, \gamma_2, \ldots, \gamma_{2g}$, and denote 
some properly embedded arcs of $P_g$ by 
$\gamma_1, \gamma_3, \ldots, \gamma_{2g+1}$, 
$\beta_4, \ldots, \beta_{2g-2}$ and 
$\beta'_4, \ldots, \beta'_{2g-2}$. 
On $\partial (P_g \times[-1,1]) = \Sigma_g$, we define 
$c_{2i-1} = \partial (\gamma_{2i-1} \times [-1,1])$ $(1 \leq i \leq g+1)$, 
$b_{2j} = \partial (\beta_{2j}\times [-1,1])$, 
$b'_{2j} = \partial (\beta'_{2j}\times [-1,1])$ $(2 \leq j \leq g-1)$, and 
$c_{2k} = \gamma_{2k} \times \{ 0 \}$ $(1 \leq k \leq g)$. 
In Figures \ref{fig:circle} and \ref{fig:b-circle}, these circles are 
illustrated and some of them are oriented. 

\begin{figure}[ht!]\small
\centering
\SetLabels
\B(0.12*1.01) $x_1$ \\ \B(0.28*1.01) $x_2$ \\ \B(0.44*1.01) $x_3$ \\
\B(0.86*1.01) $x_g$ \\
(0.18*0.19) $y_1$ \\ (0.34*0.19) $y_2$ \\ (0.50*0.19) $y_3$ \\
(0.92*0.19) $y_g$ \\
\endSetLabels
\strut\AffixLabels{\includegraphics[height=2.2cm]{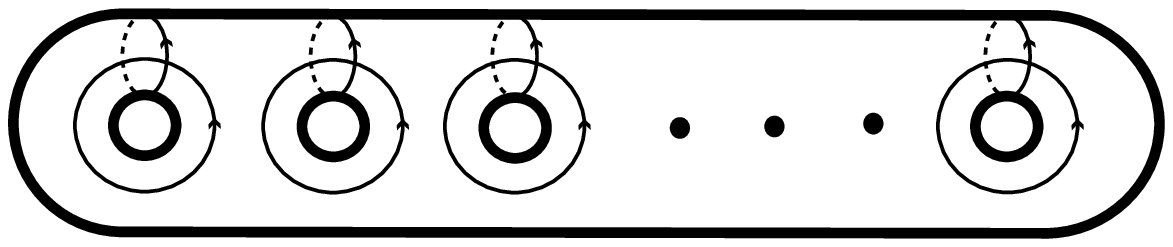}}
\nocolon
\caption{}
\label{fig:basis}
\end{figure}

We set a basis of $H_1 (\Sigma_g; \Bbb Z)$ as in Figure \ref{fig:basis}, 
where 
$x_1 = [c_1$ with opposite orientation $]$, 
$x_i = [b_{2i} \text{ with proper orientation }]$ ($2 \leq i \leq g-1$), 
$x_g = [c_{2g+1}]$, and 
$y_i = [c_{2i} \text{ with opposite orientation }]$. 

A map $q\,:\, H_1 (\Sigma_g; {\mathbb Z}_2) \to {\mathbb Z}_2$ 
is called a {\it quadratic form\/} with respect to 
the ${\Bbb Z}_2$-homology intersection form $(,)_2$ on $\Sigma_g$ 
(for short, ${\Bbb Z}_2$-{\it quadratic form\/} on $\Sigma_g$) 
if $q(x+y) = q(x) + q(y) + (x,y)_2$, for each pair of elements $x$, $y$ of 
$H_1(\Sigma_g; {\Bbb Z}_2)$. 
For the basis $\{x_1, y_1, \ldots, x_g, y_g \}$ introduced above, 
we define $Arf(q) = \sum_{i=1}^g q(x_i) q(y_i)$. 
We call a ${\Bbb Z}_2$-quadratic form $q$ {\it even\/} quadratic from 
(resp. {\it odd\/} quadratic form) if $Arf(q) = 0$ (resp. $Arf(q) = 1$). 
We define 
$$\spin{g}{q} = 
\left\{ \psi \in {\mathcal{M}}_g \; \left| \;
q (\psi_*(x)) = q(x) 
\text{ for any } x \in H_1 (\Sigma_{g}; {\Bbb Z}_2)
\right.
\right\}.   
$$
As is shown in \cite{Rourke-Sullivan}, 
for two ${\Bbb Z}_2$-quadratic forms $q$, $q'$ on $\Sigma_g$, 
if $Arf(q)$ $=$ $Arf(q')$, then there is an element $\psi'$ $\in$ ${\mathcal M}_g$ 
so that $q (\psi'_*(x)) = q'(x)$ for any $x$ $\in$ $H_1(\Sigma_g; {\Bbb Z}_2)$. 
Therefore, if $Arf(q)$ $=$ $Arf(q')$, then 
$\spin{g}{q}$ and $\spin{g}{q'}$ are conjugate in $\mathcal{M}_g$. 
By the definition of ${\Bbb Z}_2$-quadratic from, values of a quadratic form is 
completely determined by its value for the basis of $H_1(\Sigma_g; {\Bbb Z}_2)$. 
Let $q_0$ and $q_1$ be ${\Bbb Z}_2$-quadratic forms so that 
$q_0(x_i) = q_0(y_i) =0$ for $1 \leq i \leq g$, 
$q_1(x_1) = q_1(y_1) =1$ and $q_1(x_j) = q_1(y_j) = 0$ for $2 \leq j \leq g$. 
Then $q_0$ is an even quadratic form and $q_1$ an odd quadratic from. 
If $q$ is even, then $\spin{g}{q}$ is conjugate to $\spin{g}{q_0}$ in ${\mathcal M}_g$, 
on the other hand, if $q$ is odd, then $\spin{g}{q}$ is conjugate to $\spin{g}{q_1}$ 
in ${\mathcal M}_g$. 
Hence, for the sake of getting some information about groups $\spin{g}{q}$, 
it suffices to consider only on $\spin{g}{q_0}$ and $\spin{g}{q_1}$. 
The group $\spin{g}{q_0}$ is called the {\it even spin mapping class group\/}, 
and the group $\spin{g}{q_1}$ is called the {\it odd spin mapping class group\/}.
The spin mapping class group is defined by Harer \cite{Harer}, \cite{Harer2}. 
In \cite{Hirose2}, we get a system of generators for $\spin{g}{q_0}$. 
In this section, we will obtain a system of generators for $\spin{g}{q_1}$.   

Let $M$ be a simply connected smooth closed oriented 4-manifold,  
$(M, F)$ a knotted characteristic surface and 
$q_F$ the Rokhlin quadratic form for $(M,F)$. 
Rokhlin \cite{Rokhlin} showed (see also \cite{Matsumoto} and \cite{F-K}), 
$$
Arf(q_F) \equiv \frac{\sigma (M) - F \cdot F}{8} \mod 2 , 
$$ 
where $\sigma (M)$ is the signature of $M$. 
By the above formula, we can see 
$q_{K_3 \# \Sigma_{g-1}}$ is an odd quadratic form. 
Hence, we get a system of generators for 
\linebreak
$\spin{g}{q_{K_3 \# \Sigma_{g-1}}}$
from that for $\spin{g}{q_1}$. 

We introduce some notations used for describing a system of generators 
for $\spin{g}{q_1}$. 
For a simple closed curve $a$ on $\Sigma_g$, $T_a$ denotes 
the Dehn twist about $a$. 
The order of composition of maps is the functional one: 
$T_b T_a$ means we apply $T_a$ first, then $T_b$. 
For elements $a$, $b$ and $c$ of a group, we write 
$\inv{c} = c^{-1}$, and $a*b = ab \inv{a}$. 
We define some elements of $\mathcal{M}_g$ as follows:
$$
\begin{aligned}
&C_i = T_{c_i},\  B_i = T_{b_i},\  {B'}_i = T_{b'_i}, \\
&X_i = C_{i+1} C_i \inv{C_{i+1}},\  
X^*_i = \inv{C_{i+1}} C_i C_{i+1} \ \ (4 \leq i \leq 2g),\\ 
&Y_{2j} = C_{2j} B_{2j} \inv{C_{2j}},\  
Y^*_{2j} = \inv{C_{2j}} B_{2j} C_{2j} \ \ (2 \leq j \leq g-1),\\ 
&D_i = C_i^2 \ \ (1 \leq i \leq 2g+1), \\ 
&DB_{2j} = B_{2j}^2 \ \ (2 \leq j \leq g-1), \\ 
&T_1 = B_4 C_5 C_7 {\cdots} C_{2g+1}.
\end{aligned}
$$
When $g \geq 3$, $G_g$ denotes
the subgroup of $\mathcal{M}_g$ generated by $C_1$, $C_2$, $C_3$, 
$X_i$ $(4 \leq i \leq 2g)$, $Y_{2j}$ $(2 \leq j \leq g-1)$, 
$D_i$ $(1 \leq i \leq 2g+1)$, $DB_{2j}$ $(2 \leq j \leq g-1)$, and $T_1$. 
It is clear that $X^*_i$ and $Y^*_{2j}$ are elements 
of $G_g$. 
When $g=2$, the subgroup of $\mathcal{M}_2$ generated by 
$C_1$, $C_2$, $C_3$, $X_4$, and $D_j$ $(1 \leq j \leq 5)$ is denoted by $G_2$. 
For two simple closed curves $l$ and $m$ on $\Sigma_g$, 
$l$ and $m$ are called {\it $G_g$-equivalent\/} (denoted by 
$l \Ggeq m$) if there is an element $\phi$ of $G_g$ 
such that $\phi(l) = m$. 

We show that $G_g = \oddspin{g}$. That is to say, 
we show, 
\begin{thm}\label{thm:generator-Gg}
%
If $g=2$, $\oddspin{2}$ is generated by $C_1$, $C_2$, $C_3$, 
$X_4$, and $D_j$ ($1 \leq j \leq 5$). 
If $g \geq 3$, $ \oddspin{g}$ is generated by 
$C_1$, $C_2$, $C_3$, $X_i$ ($4 \leq i \leq 2g$), 
$Y_{2j}$ ($2 \leq j \leq g-1$), 
$D_k$ ($1 \leq k \leq 2g+1$), 
$DB_{2l}$ ($2 \leq l \leq g-1$), and $T_1$. \end{thm}
We prove Theorem \ref{thm:generator-Gg} by using the same method as in the proof of 
Theorem 3.1 in \cite{Hirose2}. 
By an easy calculation, we can check that each generator of $G_g$ is an element of 
$\oddspin{g}$, therefore, $G_g$ $\subset$ $\oddspin{g}$. 
Hence, we should show $\oddspin{g}$ $\subset$ $G_g$. 
In the case where $g=2$, we use the Reidemeister-Schreier method to show 
$\oddspin{g}$ $\subset$ $G_2$ (\S \ref{subsec:genus2}). 
In the case where $g \geq 3$, we use other method to show $\oddspin{g}$ $\subset$ $G_g$. 
Here, we present this method in outline. 

The integral symplectic group is denoted by $\Symp (2g, \mathbb{Z})$ and 
the $\mathbb{Z}_2$ symplectic group by $\Symp (2g, \mathbb{Z}_2)$. 
The generators of these groups are known (on $\Symp (2g, \mathbb{Z})$ see 
for example \cite{Hua-Reiner}, 
on $\Symp (2g, \mathbb{Z}_2)$ see for example \cite[Chap.3]{Grove}), 
and these generators are induced by the action of $\mathcal{M}_g$ on 
$H_1(\Sigma_g, \mathbb{Z})$ or $H_1(\Sigma_g, \mathbb{Z}_2)$. 
Therefore, the homomorphism $\Phi \co \mathcal{M}_g \to \Symp (2g, \mathbb{Z})$, 
defined by the action of $\mathcal{M}_g$ on $H_1(\Sigma_g, \mathbb{Z})$, 
is a surjection, 
and $\Psi \co \Symp (2g, \mathbb{Z}) \to \Symp (2g, \mathbb{Z}_2)$, 
defined by changing the coefficient from $\mathbb{Z}$ to $\mathbb{Z}_2$, 
is a surjection. 
In \S \ref{subsec:Step1}, we show $\ker \Phi$ $\subset$ $G_g$. 
In \S \ref{subsec:Step2}, we introduce a finite system of generators for 
$\ker \Psi$, and, for each generator, we show that one of its inverse by $\Phi$ is 
an element of $G_g$. 
Hence, we conclude $\ker \Psi \circ \Phi$ $\subset$ $G_g$. 
In \S \ref{subsec:Step3}, we introduce a finite system of generators for 
$\Psi \circ \Phi (\oddspin{g})$, and, for each generator, we show that 
one of its inverse by $\Psi \circ \Phi$ is an element of $G_g$. 
As a consequence, we show $\oddspin{g}$ $\subset$ $G_g$. 

\subsection{Step 1 for the case where $g \geq 3$}\label{subsec:Step1}  
%
%
There is a natural surjection $\Phi \co \mathcal{M}_g \to \Symp (2g, {\Bbb Z})$ 
defined by the action of $\mathcal{M}_g$ on 
$H_1(\Sigma_g; {\Bbb Z})$. The kernel of $\Phi$ is denoted by 
$\mathcal{I}_g$ and called {\it the Torelli group\/}. 
In this subsection, we prove the following lemma:
\begin{lem}\label{lem:Torelli}
%
The Torelli group $\mathcal{I}_g$ is a subgroup of 
$G_g$. 
\end{lem}
\begin{figure}[ht!]
\centering
\includegraphics[height=2cm]{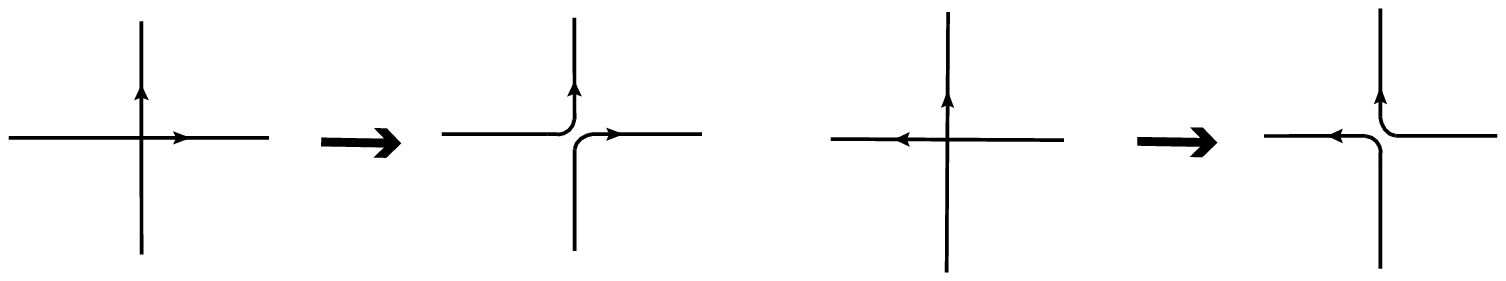}
\nocolon
\caption{}
\label{fig:smoothing}
\end{figure}
\begin{figure}[ht!]\small
\centering
\SetLabels
\B(0.24*1.01) $c_{i_2} + \cdots + c_{i_3 -1}$ \\
\B(0.83*1.01) $c_{i_r} + \cdots + c_{i_{r+1}-1}$ \\
(0.07*0.08) $c_{i_1} + \cdots + c_{i_2-1}$ \\
(0.67*0.08) $c_{i_{r-1}} + \cdots + c_{i_r -1}$ \\
\L(1.01*0.8) $b$ \\ \L(1.01*0.3) $a$ \\
\endSetLabels
\strut\AffixLabels{\includegraphics[height=4.5cm]{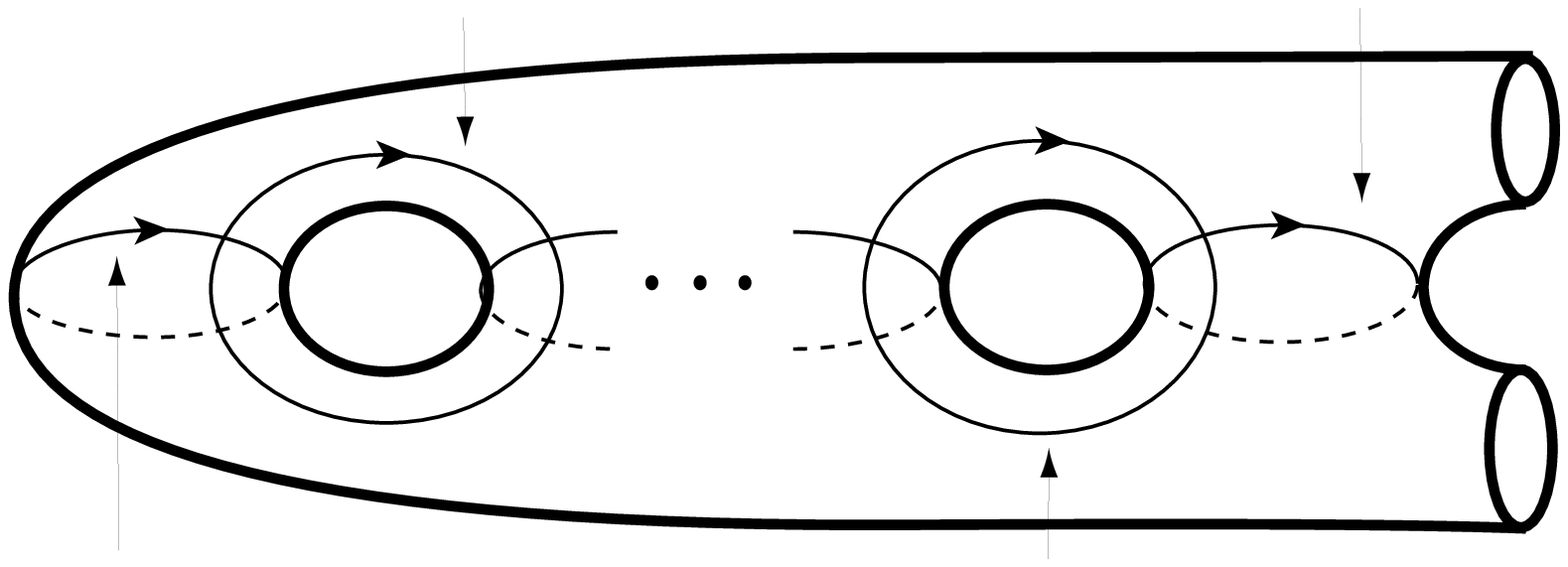}}
\nocolon
\caption{}
\label{fig:regular-nhd}
\end{figure}
Johnson \cite{Johnson1} showed that, when $g$ is larger than or equal to 
$3$, $\mathcal{I}_g$ is finitely generated. We review his result. 
For oriented simple closed curves shown in 
Figure \ref{fig:circle}, we refer to $(c_1, c_2,\ldots, c_{2g+1})$ and 
$(c_{\beta}, c_5, \ldots, c_{2g+1})$ as {\it chains\/}. 
For oriented simple closed curves $d$ and $e$ which intersect transversely 
in one point, we construct an oriented simple closed curve $d+e$ 
from $d \cup e$ as follows: 
choose a disk neighborhood of the intersection point and in it 
make a replacement as indicated in Figure \ref{fig:smoothing}. 
For a consecutive subset $\{ c_i, c_{i+1}, \ldots, c_j \}$ 
of a chain, let $c_i + \cdots + c_j$ be the oriented simple closed 
curve constructed by repeated applications of the above operations. 
Let $(i_1, \ldots, i_{r+1})$ be a subsequence of 
$(1, 2,\ldots, 2g+2)$ (resp. $(\beta, 5, \ldots, 2g+2)$). 
We construct the union of circles 
$\mathcal{C}$ $=$ $(c_{i_1}+\cdots+c_{i_2-1}) \cup (c_{i_2}+\cdots+c_{i_3-1}) 
\cup \cdots \cup (c_{i_r}+\cdots+c_{i_{r+1}-1})$. 
If $r$ is odd, a regular neighborhood of $\mathcal{C}$ is homeomorphic to 
the compact surface indicated in Figure \ref{fig:regular-nhd} whose boundaries 
are $a$ and $b$. 
Let $\phi$ $=$ $T_b T_a^{-1}$, then $\phi$ is an element of $\mathcal{I}_g$. 
We denote $\phi$ by $[i_1, \ldots, i_{r+1}]$, and call this 
{\it the odd subchain map\/} of $(c_1, c_2,\ldots, c_{2g+1})$ 
(resp. $(c_{\beta}, c_5, \ldots, c_{2g+1})$) with {\em length\/} $r+1$. 
Johnson \cite{Johnson1} showed the following theorem: 
\begin{thm}\label{thm:Johnson-gen}{\rm\cite[Main
Theorem]{Johnson1}}\qua 
%
For $g \geq 3$, the odd subchain maps of the two chains 
$(c_1, c_2,\ldots, c_{2g+1})$ and $(c_{\beta}, c_5, \ldots, c_{2g+1})$
generate $\mathcal{I}_g$. 
\end{thm}
We use the following results by Johnson \cite{Johnson1}.  
\begin{lem}\label{lem:Johnson-action}{\rm\cite{Johnson1}}\qua
%
(a) $C_j$ commutes with $[i_1, i_2, \cdots]$ if and only if $j$ and 
$j+1$ are either both contained in or are disjoint from the $i$'s. 
\newline
(b) If $i \not= j+1 $, then $\inv{C_j} * [ \cdots, j, i, \cdots ]$ 
$=$ $[ \cdots, j+1, i, \cdots ]$. 
\newline
(c) If $k \not= j$, then $C_j*[\cdots, k, j+1, \cdots]$ $=$ 
$[\cdots, k, j, \cdots]$. 
\newline
(d) $[1,2,3,4] [1,2,5,6,\ldots, 2n] B_4*[3,4,5,\ldots,2n]$ $=$ 
$[5,6,\ldots,2n][1, 2, 3, 4, \ldots,$ $2n]$, where $3 \leq n \leq g$.      
\end{lem}
\begin{rem}
Johnson showed (d) only in the case where $n=g$. 
But we can apply the proof of Lemma 10 of \cite{Johnson1} for the case where 
$3 \leq n < g$, since we can regard each surfaces in Figure 18 of \cite{Johnson1} 
as a surface of genus $n$ which is a submanifold of $\Sigma_g$. 
\end{rem}

We prove that any odd subchain map of 
$(c_1, c_2, c_3, \ldots, c_{2g+1})$ or 
$(c_{\beta}, c_5, c_6, \ldots,$ $c_{2g})$ is a product of 
elements of $G_g$. 
The following lemma shows that 
any odd subchain map of $(c_{\beta}, c_5, c_6, \ldots, c_{2g})$ 
is a product of an odd subchain map of 
$(c_1, c_2, c_3, \ldots, c_{2g+1})$ and elements of $G_g$. 
\begin{lem}\label{lem:beta-to-straight}
%
For any odd subchain map $h$ of $(c_{\beta}, c_5, c_6, \ldots, c_{2g+1})$, 
there is an element $g$ of $G_g$ such that $g*h$ is an odd subchain map of 
$(c_1, c_2, c_3, \ldots, c_{2g+1})$. 
\end{lem}
\begin{proof}
If there is not $\beta$ in the sequence which define $h$, then $h$ is 
an odd subchain map of 
$(c_1, c_2, c_3, \ldots, c_{2g+1})$. 
Hence, it suffices to treat the case where the sequence defining $h$ 
includes $\beta$. 
If $g$ $=$ $C_{2g+1}^{\epsilon_{2g+1}} \cdots C_7^{\epsilon_7} 
C_5^{\epsilon_5} B^{-1}$ ($\epsilon_i = \pm 1$), then, under 
any choice of signs of $\epsilon_i$, $g$ $\in$ $G_g$. 
We can choose signs of $\epsilon_i$ so that $g*h$ is an odd 
subchain map of $(c_1, c_2, c_3, \ldots, c_{2g+1})$. 
\end{proof}
From here to the end of this subsection, 
odd subchain maps mean only those of $(c_1, c_2, c_3, \ldots, c_{2g+1})$. 
The following lemma shows that any odd subchain map, whose length is 
at least 5 and which begins from $1,2,3,4,5$, is a product 
of shorter odd subchain maps and elements of $G_g$. 
\begin{lem}\label{lem:chain-shortening}
%
For any $6 \leq n_6 < n_7 < \cdots < n_{2k} \leq 2g+2$, 
\begin{align*}
(C_4^2)*[1,2,3,5] & [1,2,4,n_6,n_7,\ldots,n_{2k}] 
(C_4 B_4 \inv{C_4})*[3,4,5,n_6, n_7, \ldots, n_{2k}] = \\
&= [4,n_6,n_7,\ldots,n_{2k}][1,2,3,4,5,n_6,n_7,\ldots,n_{2k}]
\end{align*}
\end{lem}
\begin{proof}
By (a) of Lemma \ref{lem:Johnson-action}, 
$\inv{C_4}*[3,4,5,\ldots,2k] = [3,4,5,\ldots,2k]$, and by (d) of Lemma 
\ref{lem:Johnson-action}, 
\begin{align*}
[1,2,3,4][1,2,5,6, \ldots, 2k] &\cdot (B_4 \inv{C_4})*[3,4,5,\ldots,2k] = \\
&= [5,6, \ldots, 2k][1,2,3,4,\ldots,2k].
\end{align*}
By applying $C_4$ to the above equation and remarking that 
$C_4*[1,2,3,4]$ $=$ $(C_4^2)*(\inv{C_4}*[1,2,3,4])$ $=$ 
$(C_4^2)*[1,2,3,5]$,  we get, 
\begin{align*}
(C_4^2)*[1,2,3,5] \cdot [1,2,4,6,\cdots,2k] \cdot
&(C_4 B_4 \inv{C_4})*[3,4,5,6,\cdots,2k] = \\ 
&= [4,6,7,\cdots,2k][1,2,3,4,5,6,\cdots,2k]. 
\end{align*}
After proper applications of $\inv{C_6}$, $\inv{C_7}$, \ldots, 
$\inv{C_{2g+1}}$, we get the equation we need. 
\end{proof}
\begin{lem}\label{lem:Gg-action}
%
(1) When $i-k \geq 3$, 
$(\inv{C_{i-1}} C_{i-2} C_{i-1})*[\ldots,k,i,j,\ldots] 
=[\ldots,k,i-2,j,\ldots]$. \newline
(2) When $i-k \geq 2$, 
$(C_i C_{i-1} \inv{C_i})*[\ldots,k,i,i+1,\ldots]
=[\ldots,k,i-1,i,\ldots]$. 
\end{lem}
\begin{proof}
Lemma \ref{lem:Johnson-action} shows (1) and (2). 
\end{proof}

For any odd subchain map $[i_1, i_2, \ldots, i_r]$, we introduce a notation 
$[[\tau_1, \tau_2, \ldots,$ $\tau_{2g+2}]]$ :
$\tau_k =1$ if $k$ is a member of $\{ i_1, i_2, \ldots, i_r\}$, 
and $\tau_k = 0$ if $k$ is not a member of 
$\{ i_1, i_2, \ldots, i_r\}$. 
For $[[\tau_1, \tau_2, \ldots, \tau_{2g+2}]]$, $\tau_i$ ($1 \leq i \leq 2g+2$) 
is called the $i$-th {\em tack\/} of 
$[[\tau_1, \tau_2, \ldots, \tau_{2g+2}]]$, 
and if $\tau_i = 0$ (resp. $1$) then $\tau_i$ is called a $0$-tack (resp. a $1$-tack). 
The number of $1$-tacks in $[[\tau_1, \tau_2, \ldots, \tau_{2g+2}]]$ is called 
the {\em length\/} of $[[\tau_1, \tau_2, \ldots, \tau_{2g+2}]]$. 
Lemma \ref{lem:Gg-action} (1) means that, when $k \geq 3$, 
if there is a sequence of $0$-tacks which begins from the $k+1$-st tack and 
whose length is at least $2$, then the $1$-tack subsequent to this $0$-tack 
sequence is moved to left by $2$-steps under the action of $G_g$. 
Lemma \ref{lem:Gg-action} (2) means that, when $k \geq 3$, 
if there is a sequence of $0$-tacks which begins from the $k+1$-st tack and 
whose length is at least $1$, then the adjacent two $1$-tacks 
subsequent to this $0$-tack sequence is moved to left by $1$-step 
under the action of $G_g$. 
Therefore, for any $[[\tau_1, \tau_2, \ldots, \tau_{2g+2}]]$, 
we see, 
$$
[[\tau_1, \tau_2, \ldots, \tau_{2g+2}]] \Ggeq 
[[\tau_1, \tau_2, \tau_3, 1, \ldots, 1, 0,1, \ldots, 0,1,0, \ldots, 0]] ,
$$
where $1, \ldots, 1$ is a sequence of $1$-tacks ($b$ denotes the length of this 
sequence), $0,1,\ldots,0,1$ is a sequence arranged $0$-tacks and 
$1$-tacks alternatively ($t$ denotes the number of $1$-tacks in this sequence ), 
$0, \ldots, 0$ is a sequence of $0$-tacks. 
Since $C_1, C_2, C_3$ $\in G_g$, if there is one $1$-tack among 
$\tau_1, \tau_2, \tau_3$, then $[[\tau_1, \tau_2, \tau_3, \ldots]]$ $\Ggeq$ 
$[[1,0,0,\ldots]]$, if there are two $1$-tacks among 
$\tau_1, \tau_2, \tau_3$, then $[[\tau_1, \tau_2, \tau_3, \ldots]]$ $\Ggeq$ 
$[[1,1,0,\ldots]]$. The number of $1$-tacks in 
$\tau_1, \tau_2, \tau_3$ is denoted by $h$. 
\begin{lem}\label{lem:reduction1}
%
Any odd subchain map is a product of elements of $G_g$ and 
the odd subchain maps whose $h$ and $b$ are 
(1) $h=3, b=1$, (2) $h=3, b=0$, (3) $h=2, b=0$, (4) $h=1, b=0$, 
(5) $h=0, b=0$. 
\end{lem}
\begin{proof}
We treat the case where $h=3$. If $b \geq 2$, 
by Lemma \ref{lem:chain-shortening}, this odd subchain map is a product of 
elements of $G_g$ and shorter odd subchain maps. 

We treat the case where $h=2$. If $b \geq 3$, 
$$
[[1,1,0,1,1,1,\ldots ]] \underset{C_3}{\longrightarrow} 
[[1,1,1,0,1,1,\ldots ]] \underset{\text{Lemma \ref{lem:Gg-action}(2)}}
{\longrightarrow} [[1,1,1,1,1,0,\ldots]], $$
by Lemma \ref{lem:chain-shortening}, the last odd subchain map is a product of 
elements of $G_g$ and shorter odd subchain maps. 
If $b=2$, 
$$
[[1,1,0,1,1,0,\ldots]] \underset{C_3}{\longrightarrow} 
[[1,1,1,0,1,0,\ldots]], 
$$
the last odd subchain map is in the case where $h=3, b=0$. 
If $b=1$, $t$ should be at least $1$, and 
\begin{align*}
[[1,1,0,1,0,1,0,\ldots]] &\underset{C_3}{\longrightarrow} 
[[1,1,1,0,0,1,0,\ldots]] \\ 
&\underset{\text{Lemma \ref{lem:Gg-action}(1)}}
{\longrightarrow} [[1,1,1,1,0,0,0, \ldots]], 
\end{align*}
the last odd subchain map is in the case where $h=3, b=1$. 

We treat the case where $h=1$. If $b \geq 5$, 
\begin{align*}
&[[1,0,0,1,1,1,1,1,\ldots ]] \underset{C_2 C_3}{\longrightarrow} 
[[1,1,0,0,1,1,1,1,\ldots ]] \\ 
&\underset{\Shift}{\longrightarrow} [[1,1,0,1,1,0,1,1,\ldots ]] 
\underset{\Shift}{\longrightarrow}
[[1,1,0,1,1,1,1,0,\ldots ]] \\
&\underset{C_3}{\longrightarrow} 
[[1,1,1,0,1,1,1,0,\ldots ]] \underset{\Shift}{\longrightarrow}
[[1,1,1,1,1,0,1,0,\ldots ]], 
\end{align*}
by Lemma \ref{lem:chain-shortening}, the last odd subchain map is a 
product of elements of $G_g$ and the shorter odd subchain maps. 
If $b = 4$, 
\begin{align*}
&[[1,0,0,1,1,1,1,0,\ldots ]] \underset{C_2 C_3}{\longrightarrow}
[[1,1,0,0,1,1,1,0,\ldots ]] \\
&\underset{\Shift}{\longrightarrow}
[[1,1,0,1,1,0,1,0,\ldots ]] \underset{C_3}{\longrightarrow} [[1,1,1,0,1,0,1,0,\ldots ]], 
\end{align*}
the last odd subchain map is in the case where $h=3, b=0$. 
If $b=3$ and $t=0$, 
\begin{align*}
&[[1,0,0,1,1,1,0,0,\ldots ]] \underset{C_2 C_3}{\longrightarrow}
[[1,1,0,0,1,1,0,0,\ldots ]] \\
&\underset{\Shift}{\longrightarrow}
[[1,1,0,1,1,0,0,0,\ldots ]] 
\underset{C_3}{\longrightarrow} [[1,1,1,0,1,0,0,0,\ldots ]], 
\end{align*}
the last odd subchain map is in the case where $h=3, b=0$. 
If $b=3$ and $t \geq 2$, 
\begin{align*}
&[[1,0,0,1,1,1,0,1,\ldots ]] \underset{C_2 C_3}{\longrightarrow}
[[1,1,0,0,1,1,0,1,\ldots ]] \\ 
&\underset{\Shift}{\longrightarrow}
[[1,1,0,1,1,0,0,1,\ldots ]] 
\underset{\Jump}{\longrightarrow} [[1,1,0,1,1,1,0,0,\ldots ]] \\
&\underset{C_3}{\longrightarrow} [[1,1,1,0,1,1,0,0,\ldots ]] 
\underset{\Shift}{\longrightarrow} [[1,1,1,1,1,0,0,0,\ldots ]], 
\end{align*}
by Lemma \ref{lem:chain-shortening}, the last odd subchain map is 
a product of elements of $G_g$ and shorter odd subchain maps. 
If $b=2$, 
$$
[[1,0,0,1,1,0, \ldots ]] 
\underset{C_2 C_3}{\longrightarrow}
[[1,1,0,0,1,0,\ldots]], 
$$
the last odd subchain map is in the case where $h=2, b=0$. 
If $b=1$, $t$ should be at least $2$, 
\begin{align*}
&[[1,0,0,1,0,1,0,1,\ldots]] \underset{C_2 C_3}{\longrightarrow}
[[1,1,0,0,0,1,0,1,\ldots]] \\ 
&\underset{\Jump}{\longrightarrow}
[[1,1,0,1,0,1,0,0,\ldots]], 
\end{align*}
the last odd subchain map is in the case where $h=2, b=1$, which 
we treat before. 

We treat the case where $h=0$. If $b \geq 7$, 
\begin{align*}
&[[0,0,0,1,1,1,1,1,1,1,\ldots]] \underset{C_1 C_2 C_3}{\longrightarrow} 
[[1,0,0,0,1,1,1,1,1,1,\ldots]] \\
&\underset{\Shift}{\longrightarrow}
[[1,0,0,1,1,1,1,1,1,0,\ldots]] \underset{C_2 C_3}{\longrightarrow}
[[1,1,0,0,1,1,1,1,1,0,\ldots]] \\
&\underset{\Shift}{\longrightarrow}
[[1,1,0,1,1,1,1,0,1,0,\ldots]] \underset{C_3}{\longrightarrow}
[[1,1,1,0,1,1,1,0,1,0,\ldots]] \\
& \underset{\Shift}{\longrightarrow}
[[1,1,1,1,1,0,1,0,1,0,\ldots]], 
\end{align*}
by Lemma \ref{lem:chain-shortening}, the last odd subchain map is a product 
of $G_g$ and shorter odd subchain maps. 
If $b=6$, 
\begin{align*}
&[[0,0,0,1,1,1,1,1,1,0,\ldots]] \underset{C_1 C_2 C_3}{\longrightarrow} 
[[1,0,0,0,1,1,1,1,1,0,\ldots]] \\
&\underset{\Shift}{\longrightarrow}
[[1,0,0,1,1,1,1,0,1,0,\ldots]] \underset{C_2 C_3}{\longrightarrow}
[[1,1,0,0,1,1,1,0,1,0,\ldots]] \\
&\underset{\Shift}{\longrightarrow}
[[1,1,0,1,1,0,1,0,1,0,\ldots]] \underset{C_3}{\longrightarrow}
[[1,1,1,0,1,0,1,0,1,0,\ldots]], 
\end{align*}
the last odd subchain map is in the case where $h=3, b=0$. 
If $b=5$, $t$ should be at least $1$ and, 
{\allowdisplaybreaks
\begin{align*}
&[[0,0,0,1,1,1,1,1,0,1,\ldots]] \underset{C_1 C_2 C_3}{\longrightarrow}
[[1,0,0,0,1,1,1,1,0,1,\ldots]] \\
&\underset{\Shift}{\longrightarrow} 
[[1,0,0,1,1,1,1,0,0,1,\ldots]] \underset{\Jump}{\longrightarrow}
[[1,0,0,1,1,1,1,1,0,0,\ldots]] \\
&\underset{C_2 C_3}{\longrightarrow}
[[1,1,0,0,1,1,1,1,0,0,\ldots]] \underset{\Shift}{\longrightarrow}
[[1,1,0,1,1,1,1,0,0,0,\ldots]] \\
&\underset{C_3}{\longrightarrow} 
[[1,1,1,0,1,1,1,0,0,0,\ldots]] \underset{\Shift}{\longrightarrow}
[[1,1,1,1,1,0,1,0,0,0,\ldots]], 
\end{align*}
}
by Lemma \ref{lem:chain-shortening}, the last odd subchain map is a 
product of elements of $G_g$ and shorter odd subchain maps. 
If $b=4$, 
\begin{align*}
&[[0,0,0,1,1,1,1,0,\ldots]] \underset{C_1 C_2 C_3}{\longrightarrow}
[[1,0,0,0,1,1,1,0,\ldots]] \\
&\underset{\Shift}{\longrightarrow}
[[1,0,0,1,1,0,1,0,\ldots]] \underset{C_2 C_3}{\longrightarrow} 
[[1,1,0,0,1,0,1,0,\ldots]], 
\end{align*}
the last odd subchain map is in the case where $h=2, b=0$. 
If $b=3$ and $t=1$, 
{\allowdisplaybreaks
\begin{align*}
&[[0,0,0,1,1,1,0,1,0,\ldots]] \underset{C_1 C_2 C_3}{\longrightarrow}
[[1,0,0,0,1,1,0,1,0,\ldots]] \\ 
&\underset{\Shift}{\longrightarrow}
[[1,0,0,1,1,0,0,1,0,\ldots]] \underset{\Jump}{\longrightarrow}
[[1,0,0,1,1,1,0,0,0,\ldots]] \\
&\underset{C_2 C_3}{\longrightarrow}
[[1,1,0,0,1,1,0,0,0,\ldots]] \underset{\Shift}{\longrightarrow}
[[1,1,0,1,1,0,0,0,0,\ldots]] \\
&\underset{C_3}{\longrightarrow}
[[1,1,1,0,1,0,0,0,0,\ldots]],
\end{align*}
}
the last odd subchain map is in the case where $h=3, b=0$. 
If $b=3$ and $t \not= 1$, then $t$ should be at least $3$ and, 
{\allowdisplaybreaks
\begin{align*}
&[[0,0,0,1,1,1,0,1,0,1,0,1,\ldots]] \\
&\underset{C_1 C_2 C_3}{\longrightarrow}[[1,0,0,0,1,1,0,1,0,1,0,1,\ldots]] \\
&\underset{\Shift}{\longrightarrow} [[1,0,0,1,1,0,0,1,0,1,0,1,\ldots]] \\
&\underset{\Jump}{\longrightarrow} [[1,0,0,1,1,1,0,1,0,1,0,0,\ldots]] \\
&\underset{C_2 C_3}{\longrightarrow} [[1,1,0,0,1,1,0,1,0,1,0,0,\ldots]] \\
&\underset{\Shift}{\longrightarrow} [[1,1,0,1,1,0,0,1,0,1,0,0,\ldots]] \\
&\underset{\Jump}{\longrightarrow} [[1,1,0,1,1,1,0,1,0,0,0,0,\ldots]] \\
&\underset{C_3}{\longrightarrow} [[1,1,1,0,1,1,0,1,0,0,0,0,\ldots]] \\
&\underset{\Shift}{\longrightarrow} [[1,1,1,1,1,0,0,1,0,0,0,0,\ldots]], 
\end{align*}
}
by Lemma \ref{lem:chain-shortening}, the last odd subchain map is 
a product of elements of $G_g$ and shorter odd subchain maps. 
If $b=2$, 
$$
[[0,0,0,1,1,0,\ldots]] \underset{C_1 C_2 C_3}{\longrightarrow}
[[1,0,0,0,1,0,\ldots]], 
$$
the last odd subchain map is in the case where $h=1, b=0$. 
If $b=1$, then $t$ should be at least $3$ and, 
{\allowdisplaybreaks
\begin{align*}
&[[0,0,0,1,0,1,0,1,0,1,\ldots]] \underset{C_1 C_2 C_3}{\longrightarrow} 
[[1,0,0,0,0,1,0,1,0,1,\ldots]] \\
&\underset{\Jump}{\longrightarrow}
[[1,0,0,1,0,1,0,1,0,0,\ldots]] \underset{C_2 C_3}{\longrightarrow}
[[1,1,0,0,0,1,0,1,0,0,\ldots]] \\
&\underset{\Jump}{\longrightarrow} 
[[1,1,0,1,0,1,0,0,0,0,\ldots]] \underset{C_3}{\longrightarrow}
[[1,1,1,0,0,1,0,0,0,0,\ldots]] \\
&\underset{\Jump}{\longrightarrow} 
[[1,1,1,1,0,0,0,0,0,0,\ldots]], 
\end{align*}
}
the last odd subchain map is in the case where $h=3, b=1$. 

This Lemma follows from the above case by case arguments and 
the induction on the length ($= h + b +t$) of odd subchain maps. 
\end{proof}
\begin{lem}\label{lem:reduction2}
%
Any odd subchain maps of the $6$ cases listed in Lemma \ref{lem:reduction1} 
are products of elements of $G_g$ and odd subchain maps 
$[[1,1,1,1,0,\ldots,0]]$, $[[1,1,1,0,1,0,\ldots,0]]$, 
$[[1,1,1,0,1,0,1,0,1,0,\ldots,0]]$, 
and $[[1,1,0,0,1,0,1,0,$ 
\linebreak
$\ldots, 0]]$, 
where $0, \ldots, 0$ are sequences of $0$-tacks. 
\end{lem}
\begin{proof}
\begin{figure}[ht!]
\centering
\includegraphics[height=2cm]{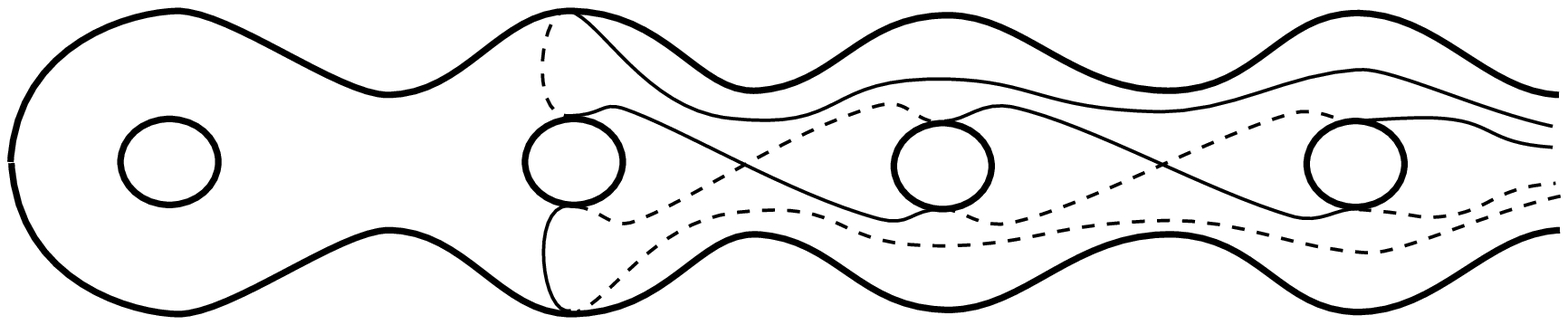}
\nocolon
\caption{}
\label{fig:commute1}
\end{figure}
\begin{figure}[ht!]
\centering
\includegraphics[height=2cm]{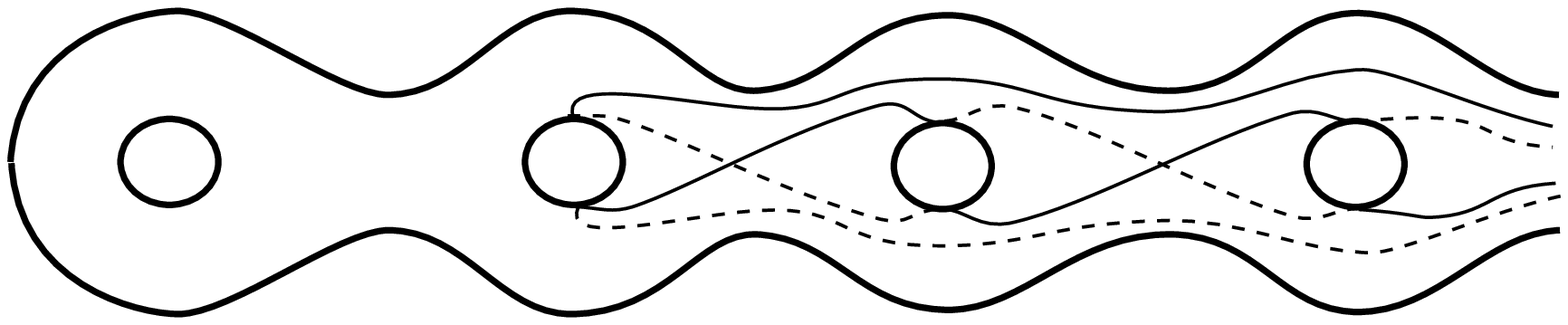}
\nocolon
\caption{}
\label{fig:commute2}
\end{figure}
By checking figures of chain maps, for examples $[[1,1,1,1,0,1,0,1,\ldots]]$ 
indicated in Figure \ref{fig:commute1} and $[[0,0,0,0,1,0,1,0,\ldots]]$ 
indicated in Figure \ref{fig:commute2}, we see that if a odd subchain map 
begins from $[[0,0,0,0,\ldots$, or $[[1,1,1,1,\ldots$, then this map commutes 
with $B_4$, hence $B_4*$ does not effect on this map. 

We treat the case where $h=3$, $b=1$. If $t=0$, then this odd subchain map is 
$[[1,1,1,1,0,\ldots,0]]$. 
If $t \not= 0$, then $t$ should be at least $2$ and, 
$$
[[1,1,1,1,0,1,0,1,\ldots]] \underset{T_1}{\longrightarrow} 
[[1,1,1,1,1,0,1,0,\ldots]], 
$$
by Lemma \ref{lem:chain-shortening}, the last odd subchain map is a product of 
elements of $G_g$ and shorter odd subchain maps. 

We treat the case where $h=3$, $b=0$. In this case, $t$ should be an odd integer 
at least $1$. 
If $t=1$, then this map is $[[1,1,1,0,1,0,\ldots,0]]$. 
If $t=3$, then this map is $[[1,1,1,0,1,0,1,0,1,0,\ldots,0]]$. 
If $t \geq 5$, 
{\allowdisplaybreaks
\begin{align*}
&[[1,1,1,0,1,0,1,0,1,0,1,0,1,0,\ldots]] \\
&\underset{\inv{C_3}}{\longrightarrow} 
[[1,1,0,1,1,0,1,0,1,0,1,0,1,0,\ldots]] \\ 
&\underset{\Shift}{\longrightarrow} [[1,1,0,0,1,1,1,0,1,0,1,0,1,0,\ldots]] \\
&\underset{\inv{C_3}\inv{C_2}}{\longrightarrow} [[1,0,0,1,1,1,1,0,1,0,1,0,1,0,\ldots]] \\
&\underset{\Shift}{\longrightarrow} [[1,0,0,0,1,1,1,1,1,0,1,0,1,0,\ldots]] \\
&\underset{\inv{C_3}\inv{C_2}\inv{C_1}}{\longrightarrow}
[[0,0,0,1,1,1,1,1,1,0,1,0,1,0,\ldots]] \\
&\underset{\Shift}{\longrightarrow} [[0,0,0,0,1,1,1,1,1,1,1,0,1,0,\ldots]] \\
&\underset{\inv{T_1}}{\longrightarrow} 
[[0,0,0,0,1,1,1,1,1,1,0,1,0,1,\ldots]] \\
&\underset{\Shift}{\longrightarrow}
[[0,0,0,1,1,1,1,1,1,0,0,1,0,1,\ldots]] \\
&\underset{C_1 C_2 C_3}{\longrightarrow}
[[1,0,0,0,1,1,1,1,1,0,0,1,0,1,\ldots]] \\
&\underset{\Jump}{\longrightarrow}
[[1,0,0,0,1,1,1,1,1,1,0,0,0,1,\ldots]] \\
&\underset{\Shift}{\longrightarrow} 
[[1,0,0,1,1,1,1,1,1,0,0,0,0,1,\ldots]] \\
&\underset{C_2 C_3}{\longrightarrow}
[[1,1,0,0,1,1,1,1,1,0,0,0,0,1,\ldots]] \\
&\underset{\Shift}{\longrightarrow}
[[1,1,0,1,1,1,1,0,1,0,0,0,0,1,\ldots]] \\
&\underset{C_3}{\longrightarrow} 
[[1,1,1,0,1,1,1,0,1,0,0,0,0,1,\ldots]] \\
&\underset{\Shift}{\longrightarrow}
[[1,1,1,1,1,0,1,0,1,0,0,0,0,1,\ldots]], 
\end{align*}
}
by Lemma \ref{lem:chain-shortening}, the last odd subchain map is a product of 
elements of $G_g$ and shorter odd subchain maps. 

We treat the case where $h=2$, $b=0$. 
In this case, $t$ should be even integer at least $2$. 
If $t=2$, this map is $[[1,1,0,0,1,0,1,0,\ldots,0]]$. 
If $t \geq 4$, 
{\allowdisplaybreaks
\begin{align*}
&[[1,1,0,0,1,0,1,0,1,0,1,0,\ldots]] \\
&\underset{\inv{C_3} \inv{C_2}}{\longrightarrow} 
[[1,0,0,1,1,0,1,0,1,0,1,0,\ldots]] \\
&\underset{\Shift}{\longrightarrow} 
[[1,0,0,0,1,1,1,0,1,0,1,0,\ldots]] \\
&\underset{\inv{C_3} \inv{C_2} \inv{C_1}}{\longrightarrow}
[[0,0,0,1,1,1,1,0,1,0,1,0,\ldots]] \\
&\underset{\Shift}{\longrightarrow} 
[[0,0,0,0,1,1,1,1,1,0,1,0,\ldots]] \\
&\underset{\inv{T_1}}{\longrightarrow} 
[[0,0,0,0,1,1,1,1,0,1,0,1,\ldots]] \\
&\underset{\Shift}{\longrightarrow} 
[[0,0,0,1,1,1,1,0,0,1,0,1,\ldots]] \\
&\underset{C_1 C_2 C_3}{\longrightarrow}
[[1,0,0,0,1,1,1,0,0,1,0,1,\ldots]] \\ 
&\underset{\Jump}{\longrightarrow}
[[1,0,0,0,1,1,1,1,0,1,0,0,\ldots]] \\
&\underset{\Shift}{\longrightarrow}
[[1,0,0,1,1,1,1,0,0,1,0,0,\ldots]] \\
&\underset{\Jump}{\longrightarrow}
[[1,0,0,1,1,1,1,1,0,0,0,0,\ldots]] \\
&\underset{C_2 C_3}{\longrightarrow}
[[1,1,0,0,1,1,1,1,0,0,0,0,\ldots]] \\
&\underset{\Shift}{\longrightarrow}
[[1,1,0,1,1,1,1,0,0,0,0,0,\ldots]] \\
&\underset{C_3}{\longrightarrow}
[[1,1,1,0,1,1,1,0,0,0,0,0,\ldots]] \\
&\underset{\Shift}{\longrightarrow}
[[1,1,1,1,1,0,1,0,0,0,0,0,\ldots]], 
\end{align*}
}
by Lemma \ref{lem:chain-shortening}, the last odd subchain map is a product of 
elements of $G_g$ and shorter odd subchain maps. 

We treat the case where $h=1$, $b=0$. 
In this case, $t$ should be an odd integer at least $3$. 
If $t=3$, 
{\allowdisplaybreaks
\begin{align*}
&[[1,0,0,0,1,0,1,0,1,0,\ldots]] \underset{\inv{C_3} \inv{C_2} \inv{C_1}}{\longrightarrow}
[[0,0,0,1,1,0,1,0,1,0,\ldots]] \\
&\underset{\Shift}{\longrightarrow} 
[[0,0,0,0,1,1,1,0,1,0,\ldots]] \underset{\inv{T_1}}{\longrightarrow} 
[[0,0,0,0,1,1,0,1,0,1,\ldots]] \\
&\underset{\Shift}{\longrightarrow} 
[[0,0,0,1,1,0,0,1,0,1,\ldots]] \underset{\Jump}{\longrightarrow}
[[0,0,0,1,1,1,0,1,0,0,\ldots]] \\
&\underset{C_1 C_2 C_3}{\longrightarrow} 
[[1,0,0,0,1,1,0,1,0,0,\ldots]] \underset{\Shift}{\longrightarrow}
[[1,0,0,1,1,0,0,1,0,0,\ldots]] \\
&\underset{\Jump}{\longrightarrow} 
[[1,0,0,1,1,1,0,0,0,0,\ldots]] \underset{C_2 C_3}{\longrightarrow}
[[1,1,0,0,1,1,0,0,0,0,\ldots]] \\
&\underset{\Shift}{\longrightarrow} 
[[1,1,0,1,1,0,0,0,0,0,\ldots]] \underset{C_3}{\longrightarrow}
[[1,1,1,0,1,0,0,0,0,0,\ldots]]. 
\end{align*}
}
If $t \geq 5$, 
{\allowdisplaybreaks
\begin{align*}
&[[1,0,0,0,1,0,1,0,1,0,1,0,\ldots]] \\
&\underset{\inv{C_3}\inv{C_2}\inv{C_1}}{\longrightarrow}
[[0,0,0,1,1,0,1,0,1,0,1,0,\ldots]] \\
&\underset{\Shift}{\longrightarrow} 
[[0,0,0,0,1,1,1,0,1,0,1,0,\ldots]] \\
&\underset{\inv{T_1}}{\longrightarrow}
[[0,0,0,0,1,1,0,1,0,1,0,1,\ldots]] \\
&\underset{\Shift}{\longrightarrow} 
[[0,0,0,1,1,0,0,1,0,1,0,1,\ldots]] \\
&\underset{\Jump}{\longrightarrow}
[[0,0,0,1,1,1,0,1,0,1,0,0,\ldots]] \\
&\underset{C_1 C_2 C_3}{\longrightarrow}
[[1,0,0,0,1,1,0,1,0,1,0,0,\ldots]] \\
&\underset{\Shift}{\longrightarrow} 
[[1,0,0,1,1,0,0,1,0,1,0,0,\ldots]] \\
&\underset{\Jump}{\longrightarrow} 
[[1,0,0,1,1,1,0,1,0,0,0,0,\ldots]] \\
&\underset{C_2 C_3}{\longrightarrow}
[[1,1,0,0,1,1,0,1,0,0,0,0,\ldots]] \\
&\underset{\Shift}{\longrightarrow} 
[[1,1,0,1,1,0,0,1,0,0,0,0,\ldots]] \\
&\underset{\Jump}{\longrightarrow} 
[[1,1,0,1,1,1,0,0,0,0,0,0,\ldots]] \\
&\underset{C_3}{\longrightarrow}
[[1,1,1,0,1,1,0,0,0,0,0,0,\ldots]] \\
&\underset{\Shift}{\longrightarrow}
[[1,1,1,1,1,0,0,0,0,0,0,0,\ldots]], 
\end{align*}
}
by Lemma \ref{lem:chain-shortening}, the last odd subchain map is a product of elements
of $G_g$ and shorter odd subchain maps. 

We treat the case where $h=0$, $b=0$. 
In this case, $t$ should be an even integer at least $4$. 
If $t=4$, 
{\allowdisplaybreaks
\begin{align*}
&[[0,0,0,0,1,0,1,0,1,0,1,0,\ldots]] \\
&\underset{\inv{T_1}}{\longrightarrow}
[[0,0,0,0,0,1,0,1,0,1,0,1,\ldots]] \\
&\underset{\Jump}{\longrightarrow} 
[[0,0,0,1,0,1,0,1,0,1,0,0,\ldots]] \\
&\underset{C_1 C_2 C_3}{\longrightarrow}
[[1,0,0,0,0,1,0,1,0,1,0,0,\ldots]] \\
&\underset{\Jump}{\longrightarrow} 
[[1,0,0,1,0,1,0,1,0,0,0,0,\ldots]] \\
&\underset{C_2 C_3}{\longrightarrow}
[[1,1,0,0,0,1,0,1,0,0,0,0,\ldots]] \\
&\underset{\Jump}{\longrightarrow} 
[[1,1,0,1,0,1,0,0,0,0,0,0,\ldots]] \\
&\underset{C_3}{\longrightarrow}
[[1,1,1,0,0,1,0,0,0,0,0,0,\ldots]] \\
&\underset{\Jump}{\longrightarrow} 
[[1,1,1,1,0,0,0,0,0,0,0,0,\ldots]]. 
\end{align*}
}
If $t \geq 6$, 
{\allowdisplaybreaks
\begin{align*}
&[[0,0,0,0,1,0,1,0,1,0,1,0,1,0,1,0,\ldots]] \\
&\underset{\inv{T_1}}{\longrightarrow} 
[[0,0,0,0,0,1,0,1,0,1,0,1,0,1,0,1,\ldots]] \\ 
&\longrightarrow \text{(as in the previous case)} \longrightarrow
[[1,1,1,1,0,0,0,0,0,0,0,0,0,1,0,1,\ldots]] \\
&\underset{\Jump}{\longrightarrow} [[1,1,1,1,0,1,0,1,\ldots]], 
\end{align*}
}
the last odd subchain map is in the case where $h=3, b=1$, which we treat before. 
\end{proof}
\begin{lem}\label{lem:formerCase}
%
The odd subchain maps $[[1,1,1,1,0,\ldots,0]]$, $[[1,1,1,0,1,0,\ldots,$ $0]]$ 
and 
$[[1,1,0,0,1,0,1,0,\ldots,0]]$ are elements of $G_g$. 
\end{lem}
\begin{proof}
In a proof of this Lemma, we use "braid relation", which is explained 
as follows. 
Let $a$ and $b$ are simple closed curves on $\Sigma_g$ intersecting transversely 
in one point, then $T_a T_b T_a^{-1}$ $=$ $T_b^{-1} T_a T_b$, in other word, 
$T_a * T_b$ $=$ $\inv{T_b}*T_a$. 

Let $b_4'$ be the simple closed curve on $\Sigma_g$ indicated 
in Figure \ref{fig:b-circle} and let $B_4' = T_{b_4'}$. 
The odd subchain map $[[1,1,1,1,0,\ldots,0]]$ is equal to $B_4 \inv{B_4'}$. 
Since $b_4'$ $=$ $C_4 C_3 C_2 C_1 C_1 C_2 C_3 C_4 (b_4)$, 
{\allowdisplaybreaks
\begin{align*}
B_4 \inv{B_4'} =& B_4 C_4 C_3 C_2 C_1 C_1 C_2 C_3 C_4 \inv{B_4} \inv{C_4} \inv{C_3} 
\inv{C_2} \inv{C_1} \inv{C_1} \inv{C_2} \inv{C_3} \inv{C_4} \\
=&(B_4 C_4 C_3 C_2)*(C_1 C_1) \cdot (B_4 C_4 C_3)*(C_2 C_2) \cdot (B_4 C_4)*(C_3 C_3) 
\cdot \\
&\cdot B_4*(C_4 C_4) 
\cdot (\inv{C_4} \inv{C_3} \inv{C_2})*(\inv{C_1} \inv{C_1}) 
\cdot (\inv{C_4} \inv{C_3})*(\inv{C_2} \inv{C_2}) \cdot \\
&\cdot \inv{C_4}*(\inv{C_3} \inv{C_3}) 
\cdot \inv{C_4} \inv{C_4}.  
\end{align*}
}
This equation means that $B_4 \inv{B_4'}$ is a product of squares Dehn twists. 
By using braid relations of ${\mathcal M}_g$, we can see that these squares of 
Dehn twists are elements of $G_g$ as follows, 
{\allowdisplaybreaks
\begin{align*}
 (B_4 C_4 C_3 C_2)*(C_1 C_1) &= 
(\inv{C_1} \cdot \inv{C_2} \cdot \inv{C_3} \cdot B_4)*(C_4 C_4) \\
&= (\inv{C_1} \cdot \inv{C_2} \cdot \inv{C_3})*
(B_4 C_4 \inv{B_4} \cdot B_4 C_4 \inv{B_4}), \\
(B_4 C_4 C_3)*(C_2 C_2) &= 
(\inv{C_2} \cdot \inv{C_3} \cdot B_4)*(C_4 C_4) \\
&= (\inv{C_2} \cdot \inv{C_3})*(B_4 C_4 \inv{B_4} \cdot B_4 C_4 \inv{B_4}), \\
(B_4 C_4)*(C_3 C_3) &= (\inv{C_3} \cdot B_4)*(C_4 C_4) = 
\inv{C_3}*(B_4 C_4 \inv{B_4} \cdot B_4 C_4 \inv{B_4}), \\  
B_4*(C_4 C_4) &= B_4 C_4 \inv{B_4} \cdot B_4 C_4 \inv{B_4}, \\
(\inv{C_4} \inv{C_3} \inv{C_2})*(C_1 C_1) &= (C_1 \cdot C_2 \cdot C_3)* (C_4 C_4), \\
(\inv{C_4} \inv{C_3})*(C_2 C_2) &= (C_2 \cdot C_3)*(C_4 C_4), \\
\inv{C_4}*(C_3 C_3) &= C_3 * (C_4 C_4). 
\end{align*}
}
Since $\inv{C_4}*[[1,1,1,1,0, \ldots, 0]] = [[1,1,1,0,1,0,\ldots,0]]$, 
{\allowdisplaybreaks
\begin{align*}
&[[1,1,1,0,1,0,\ldots,0]] = \inv{C_4}*(B_4 \inv{B_4'}) \\
&=(\inv{C_4} B_4 C_4 C_3 C_2)*(C_1 C_1) \cdot (\inv{C_4} B_4 C_4 C_3)*(C_2 C_2) 
\cdot (\inv{C_4} B_4 C_4)*(C_3 C_3) \cdot \\ 
&\cdot (\inv{C_4} B_4)*(C_4 C_4) 
\cdot (\inv{C_4} \inv{C_4} \inv{C_3} \inv{C_2})*(\inv{C_1} \inv{C_1}) 
\cdot (\inv{C_4} \inv{C_4} \inv{C_3})*(\inv{C_2} \inv{C_2}) \cdot \\
&\cdot (\inv{C_4} \inv{C_4})*(\inv{C_3} \inv{C_3}) 
\cdot \inv{C_4} * (\inv{C_4} \inv{C_4}) \\
&= (\inv{C_4} B_4 C_4 \cdot C_3 \cdot C_2)*(C_1 C_1) 
\cdot (\inv{C_4} B_4 C_4 \cdot C_3)*(C_2 C_2) \cdot \\
&\cdot (\inv{C_4} B_4 C_4)*(C_3 C_3) 
\cdot (\inv{C_4} B_4 C_4)*(C_4 C_4)
\cdot (\inv{C_4} \inv{C_4} \cdot \inv{C_3} \cdot \inv{C_2})*(\inv{C_1} \inv{C_1}) 
\cdot \\
&\cdot (\inv{C_4} \inv{C_4} \cdot \inv{C_3})*(\inv{C_2} \inv{C_2}) 
\cdot (\inv{C_4} \inv{C_4})*(\inv{C_3} \inv{C_3}) 
\cdot (\inv{C_4} \inv{C_4})
\end{align*}
}
This equation shows that $[[1,1,1,0,1,0,\ldots,0]]$ $\in$ $G_g$. 

Since $\inv{C_4} \inv{C_3} \inv{C_6} \inv{C_5} \inv{C_4} * [[1,1,1,1,0,\ldots,0]]$ $=$ 
$[[1,1,0,0,1,0,1,0,\ldots,0]]$, 
{\allowdisplaybreaks
\begin{align*}
&[[1,1,0,0,1,0,1,0,\ldots,0]] = 
\inv{C_4} \inv{C_3} \inv{C_6} \inv{C_5} \inv{C_4} *(B_4 \inv{B_4'})  \\
&=(\inv{C_4} \inv{C_3} \inv{C_6} \inv{C_5} \inv{C_4} B_4 C_4 C_3 C_2)*(C_1 C_1) \cdot 
(\inv{C_4} \inv{C_3} \inv{C_6} \inv{C_5} \inv{C_4} B_4 C_4 C_3)*(C_2 C_2) \cdot \\
&\cdot (\inv{C_4} \inv{C_3} \inv{C_6} \inv{C_5} \inv{C_4} B_4 C_4)*(C_3 C_3) \cdot 
(\inv{C_4} \inv{C_3} \inv{C_6} \inv{C_5} \inv{C_4} B_4)*(C_4 C_4) \cdot \\
&\cdot (\inv{C_4} \inv{C_3} \inv{C_6} \inv{C_5} \inv{C_4} \inv{C_4} \inv{C_3} \inv{C_2})
*(\inv{C_1} \inv{C_1}) \cdot 
(\inv{C_4} \inv{C_3} \inv{C_6} \inv{C_5} \inv{C_4} \inv{C_4} \inv{C_3})
*(\inv{C_2} \inv{C_2}) \cdot \\
&\cdot (\inv{C_4} \inv{C_3} \inv{C_6} \inv{C_5} \inv{C_4} \inv{C_4})*
(\inv{C_3} \inv{C_3}) \cdot
(\inv{C_4} \inv{C_3} \inv{C_6} \inv{C_5} \inv{C_4})* (\inv{C_4} \inv{C_4}). 
\end{align*}
}
This equation describes $[[1,1,0,0,1,0,1,0,\ldots,0]]$ as a product of squares of 
Dehn twists. 
By using braid relations of ${\mathcal M}_g$, we show that these squares of 
Dehn twists are elements of $G_g$ as follows, 
{\allowdisplaybreaks
\begin{align*}
(\inv{C_4} \inv{C_3} \inv{C_6} \inv{C_5} \inv{C_4} B_4 C_4 C_3 C_2)&*(C_1 C_1) \\
&=(\inv{C_1} \cdot \inv{C_4} B_4 C_4 \cdot \inv{C_2} \cdot C_3 \cdot \inv{C_6}C_5 C_6) * 
(C_4 C_4), \\
(\inv{C_4} \inv{C_3} \inv{C_6} \inv{C_5} \inv{C_4} B_4 C_4 C_3)&*(C_2 C_2) \\
&=(\inv{C_4} B_4 C_4 \cdot \inv{C_2} \cdot C_3 \cdot \inv{C_6}C_5 C_6) * 
(C_4 C_4), \\
(\inv{C_4} \inv{C_3} \inv{C_6} \inv{C_5} \inv{C_4} B_4 C_4)&*(C_3 C_3)
=(\inv{C_4} B_4 C_4 \cdot C_3 \cdot \inv{C_4} \inv{C_4} \cdot \inv{C_6} ) * (C_5 C_5) \\
&=(\inv{C_4} B_4 C_4 \cdot C_3 \cdot \inv{C_4} \inv{C_4}) * 
(\inv{C_6} C_5 C_6 \cdot \inv{C_6} C_5 C_6), \\
(\inv{C_4} \inv{C_3} \inv{C_6} \inv{C_5} \inv{C_4} B_4)&*(C_4 C_4) \\
&=(C_3 \cdot \inv{C_4} B_4 C_4 \cdot \inv{C_6} C_5 C_6 \cdot \inv{C_4} \inv{C_4})*
(C_3 C_3), \\
(\inv{C_4} \inv{C_3} \inv{C_6} \inv{C_5} \inv{C_4} \inv{C_4} \inv{C_3} \inv{C_2})
&*(\inv{C_1} \inv{C_1}) \\
=(C_1 \cdot & C_3 \cdot C_2 \cdot \inv{C_4} \inv{C_4} \cdot \inv{C_6} \inv{C_5} C_6 \cdot 
\inv{C_4} \inv{C_4})*(C_3 C_3), \\
(\inv{C_4} \inv{C_3} \inv{C_6} \inv{C_5} \inv{C_4} \inv{C_4} \inv{C_3})
&*(\inv{C_2} \inv{C_2}) \\
&=(C_3 \cdot C_2 \cdot \inv{C_4} \inv{C_4} \cdot \inv{C_6} \inv{C_5} C_6 \cdot 
\inv{C_4} \inv{C_4})*(C_3 C_3), \\
(\inv{C_4} \inv{C_3} \inv{C_6} \inv{C_5} \inv{C_4} \inv{C_4})&*(\inv{C_3} \inv{C_3}) 
=(C_3 \cdot \inv{C_6} C_5 C_6 \cdot \inv{C_6} C_5 C_6)*(C_4 C_4), \\
(\inv{C_4} \inv{C_3} \inv{C_6} \inv{C_5} \inv{C_4})&*(\inv{C_4} \inv{C_4})
=(C_3 \cdot \inv{C_4} C_5 C_4)* (C_6 C_6). 
\end{align*}
}
\end{proof}
\begin{lem}\label{lem:lastCase}
%
The odd subchain map $[[1,1,1,0,1,0,1,0,1,0,\ldots,0]]$ 
is an element of $G_g$. 
\end{lem}
\begin{proof}
\begin{figure}[ht!]
\centering
\includegraphics[height=2cm]{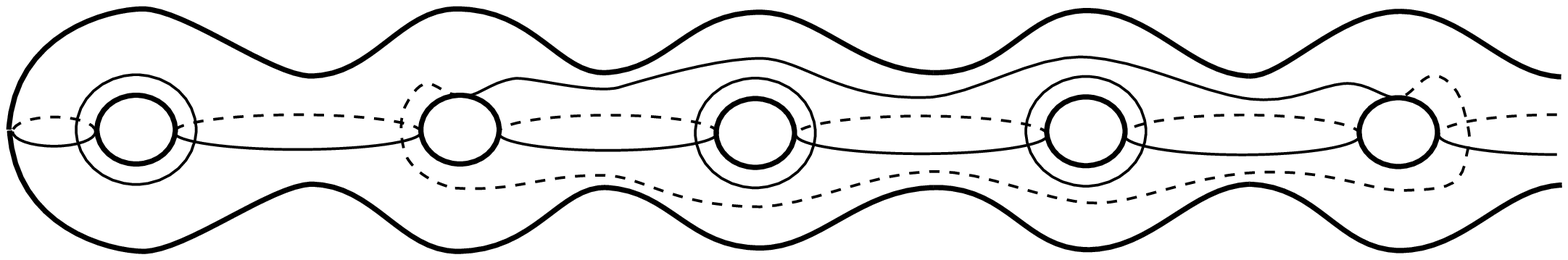}
\nocolon
\caption{}
\label{fig:complement}
\end{figure}
We can show that this odd subchain map is $G_g$-equivalent to 
$[[0,0,0,0,$ $1,1,1,1,1,1,0,\ldots,0]]$ as follows, 
{\allowdisplaybreaks
\begin{align*}
[[1,1,1,0,1,0,1,0,1,0,&\ldots,0]] \underset{\inv{C_3}}{\longrightarrow} 
[[1,1,0,1,1,0,1,0,1,0,\ldots,0]] \\
&\underset{\Shift}{\longrightarrow} 
[[1,1,0,0,1,1,1,0,1,0,\ldots,0]] \\
&\underset{\inv{C_3}\inv{C_2}}{\longrightarrow}
[[1,0,0,1,1,1,1,0,1,0,\ldots,0]] \\ 
&\underset{\Shift}{\longrightarrow} 
[[1,0,0,0,1,1,1,1,1,0,\ldots,0]] \\
&\underset{\inv{C_3}\inv{C_2}\inv{C_1}}{\longrightarrow}
[[0,0,0,1,1,1,1,1,1,0,\ldots,0]] \\
&\underset{\Shift}{\longrightarrow} 
[[0,0,0,0,1,1,1,1,1,1,0,\ldots,0]]. 
\end{align*}
}
If $g=4$, $[[0,0,0,0,1,1,1,1,1,1]] = B_4 \inv{B_4'} = [[1,1,1,1,0,0,0,0,0,0]]$, 
which we have already treated in Lemma \ref{lem:formerCase}. 
If $g \geq 5$, as we see in Figure \ref{fig:complement}, 
$$
[[0,0,0,0,1,1,1,1,1,1,0,\ldots,0]] = 
[[1,1,1,1,0,0,0,0,0,0,1,\ldots,1]], 
$$
in the notation of the last odd subchain map, $\ldots$ is a sequence of $1$-tacks. 
By Lemma \ref{lem:Gg-action} (2), 
$$
[[1,1,1,1,0,0,0,0,0,0,1,\ldots,1]] \Ggeq [[1,1,1,1,1,\ldots,1,0,0,0,0,0,0]], 
$$
which is a product of elements of $G_g$ and shorter odd subchain maps. 
\end{proof}

Therefore, Lemma \ref{lem:Torelli} is proved. 
\subsection{Step 2 for the case where $g \geq 3$}\label{subsec:Step2}
%
%
Let $\Phi_2$ be the natural homomorphism from $\mathcal{M}_g$ to 
$\Symp (2g, {\Bbb Z}_2)$ defined by the action of $\mathcal{M}_g$ 
on the ${\Bbb Z}_2$-coefficient first homology group 
$H_1(\Sigma_g ; {\Bbb Z}_2)$. 
In this section, we will show the following lemma. 

\begin{lem}\label{lem:level2-congruence}
%
$\ker \Phi_2$ is a subgroup of $G_g$. 
\end{lem}
We denote the kernel of the natural homomorphism from 
$\Symp(2g,{\Bbb Z})$ 
\linebreak
to $\Symp(2g,{\Bbb Z}_2)$ by $\Symp^{(2)} (2g)$. 
We set a basis of $H_1(\Sigma_g; {\Bbb Z})$ as in Figure \ref{fig:basis}, 
and define the intersection form $(,)$ on $H_1(\Sigma_g; {\Bbb Z})$ 
to satisfy $(x_i, y_j)=\delta_{i,j}$, $(x_i, x_j) = (y_i, y_j) = 0$ 
$(1 \leq i, j, \leq g)$. 
An element $a$ of $H_1(\Sigma_g; {\Bbb Z})$ is called 
{\it primitive\/} if there is no element $n (\not=0, \pm 1)$ of 
${\Bbb Z}$, and no element $b$ of $H_1(\Sigma_g; {\Bbb Z})$ such that 
$a = nb$. 
For a primitive element $a$ of $H_1(\Sigma_g; {\Bbb Z})$, we define 
an isomorphism $T_a: H_1(\Sigma_g; {\Bbb Z}) \to H_1(\Sigma_g; {\Bbb Z})$ 
by $T_a(v) = v + (a,v)a$. 
This isomorphism is the action of Dehn twist about a simple
closed curve  representing $a$ on $H_1(\Sigma_g; {\Bbb Z})$.  
We call $T_a^2$ {\it the square transvection about $a$\/}. 
Johnson \cite{Johnson2} showed the following result. 
\begin{lem}\label{lem:double-transvection}
%
$\Symp^{(2)} (2g)$ is generated by square transvections. 
\end{lem}
In \cite{Hirose2}, we showed, 
\begin{lem}\label{lem:level2-congruence-generator}
%
$\Symp^{(2)} (2g)$ is generated by the square transvections about 
the primitive elements 
$\sum_{i=1}^g (\epsilon_i x_i + \delta_i y_i)$, 
where $\epsilon_i =0,1$ and $\delta_i =0,1$. 
\end{lem}
\begin{figure}[ht!]\small
\centering
\SetLabels
(0.05*0.02) (-) \\ (0.2*0.02) (0) \\ (0.4*0.02) (1) \\ (0.6*0.02) (2) \\ 
(0.79*0.02) (3) \\ (0.95*0.02) (+) \\
\endSetLabels
\strut\AffixLabels{\includegraphics[height=2cm]{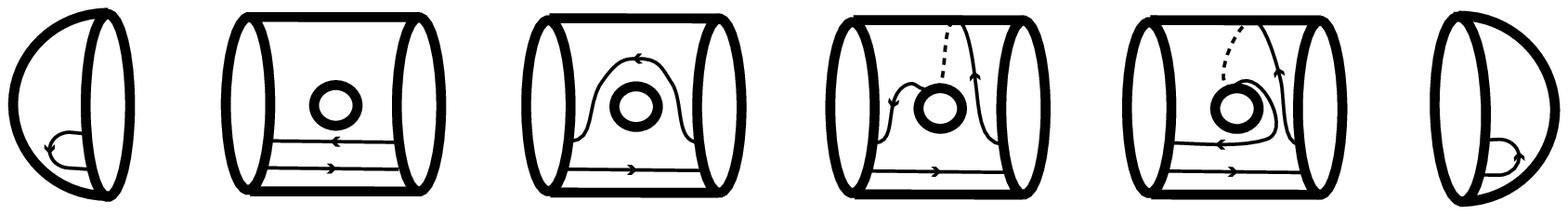}}
\nocolon
\caption{}
\label{fig:realization}
\end{figure}
For each element $[(\epsilon_1, \delta_1), \cdots, (\epsilon_g, \delta_g)]$ 
$= \sum_{i=1}^g (\epsilon_i x_i + \delta_i y_i)$ 
(where $\epsilon_i = 0,1$, $\delta_i = 0,1$) of 
$H_1(\Sigma_g; {\Bbb Z})$, we construct an oriented simple closed curve 
on $\Sigma_g$ which represent this homology class. 
For each $i$-th block, 
if $(\epsilon_i, \delta_i) = (0,0)$, we prepare (0) 
of Figure \ref{fig:realization}, 
if $(\epsilon_i, \delta_i) = (0,1)$, we prepare (1) 
of Figure \ref{fig:realization}, 
if $(\epsilon_i, \delta_i) = (1,1)$, we prepare (2) 
of Figure \ref{fig:realization}, 
if $(\epsilon_i, \delta_i) = (1,0)$, we prepare (3) 
of Figure \ref{fig:realization}. 
After that, we glue them along the boundaries and 
cap the left boundary component by (-) of Figure \ref{fig:realization} 
and the right boundary component by (+) of Figure \ref{fig:realization}. 
We denote this oriented simple closed curve on $\Sigma_g$ by 
$\{ (\epsilon_1, \delta_1), \cdots, (\epsilon_g, \delta_g) \}$. 
Here, we remark that 
the action of 
$T_{\{ (\epsilon_1, \delta_1), \cdots, (\epsilon_g, \delta_g)
\}}$ on $H_1(\Sigma_g; {\Bbb Z})$ equals  
$T_{[(\epsilon_1, \delta_1), \cdots, (\epsilon_g, \delta_g)]}$, 
and, for any $\phi$ of $\mathcal{M}_g$, 
$\phi \circ T_{\{ (\epsilon_1, \delta_1), \cdots, (\epsilon_g, \delta_g)
\}} \circ \phi^{-1}$ $=$ 
$T_{\phi (\{ (\epsilon_1, \delta_1), \cdots, (\epsilon_g, \delta_g)
\})}$. 
\begin{lem}\label{lem:level2-Gg-action}
%
For any $\{ (\epsilon_1, \delta_1), \cdots, (\epsilon_g, \delta_g) \}$,     
there is an element $\phi$ of $G_g$ such that 
%
\begin{align*}
\phi(\{ (\epsilon_1, \delta_1), \cdots, (\epsilon_g, \delta_g)\}) 
&= \{(0,0),(0,1),(0,0),(0,0), \cdots, (0,0)\} \\
\text{ or } &= \{(0,0),(1,1),(0,0),(0,0),\cdots,(0,0)\} \\
\text{ or } &= \{(0,0),(0,0),(1,1),(0,0),\cdots,(0,0)\} \\
\text{ or } &= \{(0,1),(0,0),(0,0), \cdots, (0,0)\} \\
\text{ or } &= \{(1,1),(0,0),(0,0),\cdots,(0,0)\} \\
\text{ or } &= \{(0,0),(0,0),(0,0),\cdots,(0,0)\}. 
\end{align*}
%
\end{lem}
\begin{figure}[ht!]\small
\centering
\SetLabels
\R(-0.01*0.94) (a) \\ \R(-0.01*0.75) (b) \\ \R(-0.01*0.56) (c) \\
\R(-0.01*0.35) (d) \\ \R(-0.01*0.13) (e) \\
(0.32*0.88) $\inv{X_{2i+1}^*}$ \\ (0.65*0.88) $\inv{X_{2i}^*}$ \\ 
(0.32*0.70) $\inv{X_{2i}^*}$ \\ (0.65*0.70) $\inv{D_{2i}}$ \\
(0.23*0.465) $X_{2i}$ \\ (0.49*0.465) $\inv{DB_{2i+2}}$ \\ 
(0.75*0.465) $\inv{X_{2i}^*} \inv{X_{2i+1}^*}$ \\
(0.23*0.25) $\inv{Y_{2i+2}^*}$ \\ (0.49*0.25) $\inv{X_{2i}}$ \\
(0.75*0.25) $Y_{2i}^*$ \\
(0.23*0.03) $\inv{DB_{2i+2}}$ \\ (0.49*0.03) $\inv{X_{2i+1}}$ \\
(0.75*0.03) $\inv{DB_{2i}}$ \\
\endSetLabels
\strut\AffixLabels{\includegraphics[height=10.5cm]{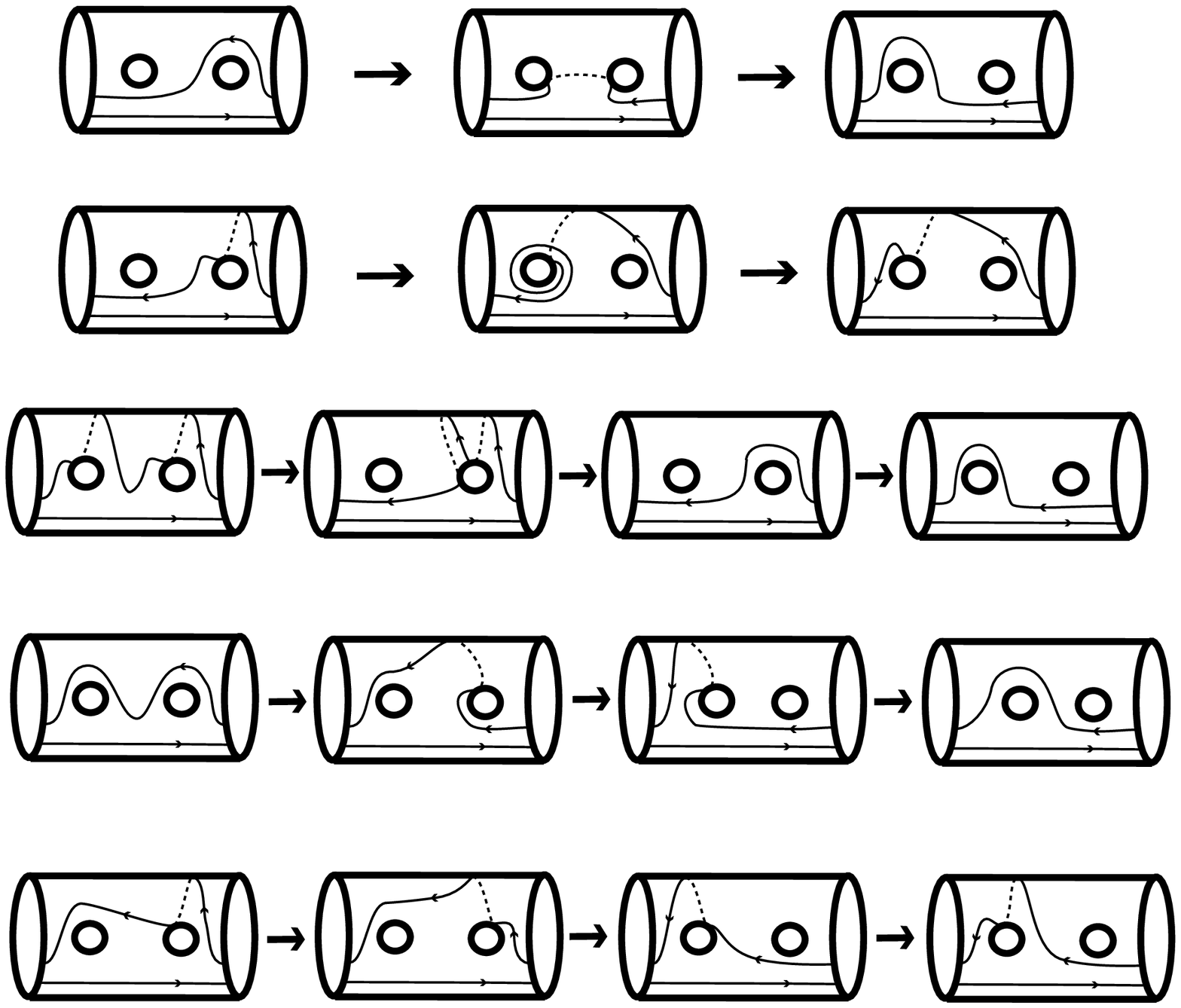}}
\nocolon
\caption{}
\label{fig:elementary}
\end{figure}
\begin{proof}
If the $i$-th block is (3), by the action of $\inv{Y_{2i}}$ if $2 \leq i \leq g-1$, 
$C_2 \inv{C_1} \inv{C_2}$ if $i=1$, and $C_{2g} \inv{C_{2g+1}} \inv{C_{2g}}$ 
if $i=g$, this block is changed to (1). 
Therefore, it suffices to show this lemma in the case where each block 
is not (3). 
First we investigate actions of elements of $G_g$ on adjacent 
blocks, say the $i$-th block and the $i+1$-st block, where $i \geq 2$. 
Each picture of Figure \ref{fig:elementary} shows the action of $G_g$ on this 
adjacent blocks. 
%
\begin{align*}
&(a) \text{ shows } \{ \Cbullet, (0,0),(0,1), \Cbullet \} \Ggeq 
\{ \Cbullet, (0,1), (0,0), \Cbullet \}, \\ 
&(b) \text{ shows } \{ \Cbullet, (0,0),(1,1), \Cbullet \} \Ggeq 
\{ \Cbullet, (1,1), (0,1), \Cbullet \}, \\
&(c) \text{ shows } \{ \Cbullet, (1,1),(1,1), \Cbullet \} \Ggeq 
\{ \Cbullet, (0,1), (0,0), \Cbullet \}, \\
&(d) \text{ shows } \{ \Cbullet, (0,1),(0,1), \Cbullet \} \Ggeq 
\{ \Cbullet, (0,1), (0,0), \Cbullet \}, \\
&(e) \text{ shows } \{ \Cbullet, (0,1),(1,1), \Cbullet \} \Ggeq 
\{ \Cbullet, (1,1), (0,0), \Cbullet \}, 
\end{align*}
where $\Cbullet$ indicates the part which is not changed by 
the action of $G_g$. 
Let $x=$ $\{ (\epsilon_1, \delta_1), \cdots, (\epsilon_g, \delta_g)\}$, 
each of whose block is $(0,0)$ or $(0,1)$ or $(1,1)$.  
If there are the $j$-th blocks $(1,1)$ $(j \geq 2)$, by (b) and (e), 
they are gathered to a sequence of $(1,1)$ blocks 
which begins from the second block. 
If there are the $j$-th blocks $(0,1)$ $(j \geq 2)$, by (a), 
they are gathered to a sequence of $(0,1)$ blocks 
subsequent to the previous sequence of $(1,1)$ blocks. 
Hence, we showed, 
$$
x \Ggeq \{(\epsilon_1, \delta_1), (1,1), \cdots, (1,1), (0,1), \cdots, 
(0,1), (0,0), \cdots, (0,0) \}. 
$$
By (a) and (d), the sequence of $(0,1)$ blocks is altered to 
$(0,1),(0,0), \cdots, (0,0)$ or $(0,0), \cdots, (0,0)$. 
By (c), the sequence of $(1,1)$ blocks is altered to 
$(1,1)$, $(0,1), (0,0), \cdots, (0,1), (0,0)$ (when the length of the sequence 
is odd) or to 
$(0,1), (0,0), \cdots, (0,1), (0,0)$ (when the length of the sequence 
is even). 
By (a) and (d), $(1,1), (0,1), (0,0), \cdots, (0,1), (0,0)$ is altered to 
$(1,1), (0,1), (0,0)$, $\cdots$, $(0,0), (0,0)$ or 
$(1,1), (0,0), (0,0), \cdots, (0,0), (0,0)$, and 
$(0,1), (0,0), \cdots, (0,1)$, $(0,0)$ to 
$(0,1), (0,0), \cdots, (0,0), (0,0)$ or 
$(0,0), (0,0), \cdots, (0,0), (0,0)$. 
Therefore, we showed, 
{\allowdisplaybreaks
\begin{align*}
x &\Ggeq 
\{(\epsilon_1, \delta_1), (1,1), (0,0), (0,0), \cdots, (0,0),
(0,0), (0,0), \cdots, (0,0) \}, \\
&\text{ or } \Ggeq 
\{(\epsilon_1, \delta_1), (1,1), (0,0), (0,0), \cdots, (0,0), 
(0,1), (0,0), \cdots, (0,0) \}, \\
&\text{ or } \Ggeq 
\{(\epsilon_1, \delta_1), (1,1), (0,1), (0,0), \cdots, (0,0), 
(0,0), (0,0), \cdots, (0,0) \}, \\
&\text{ or } \Ggeq 
\{(\epsilon_1, \delta_1), (1,1), (0,1), (0,0), \cdots, (0,0), 
(0,1), (0,0), \cdots, (0,0) \}, \\
&\text{ or } \Ggeq 
\{(\epsilon_1, \delta_1), (0,1), (0,0), (0,0), \cdots, (0,0),
(0,0), (0,0), \cdots, (0,0) \}, \\
&\text{ or } \Ggeq 
\{(\epsilon_1, \delta_1), (0,1), (0,0), (0,0), \cdots, (0,0), 
(0,1), (0,0), \cdots, (0,0) \}, \\
&\text{ or } \Ggeq 
\{(\epsilon_1, \delta_1), (0,0), (0,0), (0,0), \cdots, (0,0),
(0,0), (0,0), \cdots, (0,0) \}, \\
&\text{ or } \Ggeq 
\{(\epsilon_1, \delta_1), (0,0), (0,0), (0,0), \cdots, (0,0), 
(0,1), (0,0), \cdots, (0,0) \}. 
\end{align*}
}
In the second case, 
\begin{align*}
& \{(\epsilon_1, \delta_1),(1,1),(0,0), \cdots, (0,0), (0,1), (0,0), 
\cdots, (0,0) \} \\
& \Ggeq \{(\epsilon_1, \delta_1), (1,1), (0,1), (0,0), \cdots, (0,0) \} 
(\text{ by } (a)).  
\end{align*}
In the 4-th case, 
\begin{align*}
& \{(\epsilon_1, \delta_1),(1,1),(0,1), (0,0), \cdots, (0,0), (0,1), (0,0), 
\cdots, (0,0) \} \\
& \Ggeq \{(\epsilon_1, \delta_1), (1,1), (0,1),(0,1), (0,0), \cdots, (0,0) \} 
(\text{ by } (a)) \\
& \Ggeq \{(\epsilon_1, \delta_1), (1,1), (0,1),(0,0), (0,0), \cdots, (0,0) \} 
(\text{ by } (d)). 
\end{align*}
In the 6-th case, 
\begin{align*}
& \{(\epsilon_1, \delta_1),(0,1), (0,0), \cdots, (0,0), (0,1), (0,0), 
\cdots, (0,0) \} \\
& \Ggeq \{(\epsilon_1, \delta_1), (0,1),(0,1), (0,0), \cdots, (0,0) \} 
(\text{ by } (a)) \\
& \Ggeq \{(\epsilon_1, \delta_1), (0,1),(0,0), (0,0), \cdots, (0,0) \} 
(\text{ by } (d)). 
\end{align*}
In the 8-th case, 
\begin{align*}
& \{(\epsilon_1, \delta_1),(0,0), \cdots, (0,0), (0,1), (0,0), 
\cdots, (0,0) \} \\
& \Ggeq \{(\epsilon_1, \delta_1),  (0,1), (0,0), \cdots, (0,0) \} 
(\text{ by } (a)).  
\end{align*}
Therefore, 
\begin{align*}
x &\Ggeq 
\{(\epsilon_1, \delta_1), (1,1), (0,0), (0,0), \cdots, (0,0)\}, \\
&\text{ or } \Ggeq 
\{(\epsilon_1, \delta_1), (1,1), (0,1), (0,0), \cdots, (0,0) \}, \\
&\text{ or } \Ggeq 
\{(\epsilon_1, \delta_1), (0,1), (0,0), (0,0), \cdots, (0,0) \}, \\
&\text{ or } \Ggeq 
\{(\epsilon_1, \delta_1), (0,0), (0,0), (0,0), \cdots, (0,0) \}. 
\end{align*}
There are 7 cases remained to consider, 
{\allowdisplaybreaks
\begin{align*}
&\{(0,0), (1,1), (0,1), (0,0), \cdots, (0,0) \}, \quad 
\{(0,1), (1,1), (0,0), (0,0), \cdots, (0,0) \}, \\
&\{(0,1), (1,1), (0,1), (0,0), \cdots, (0,0) \}, \quad 
\{(0,1), (0,1), (0,0), (0,0), \cdots, (0,0) \}, \\
&\{(1,1), (1,1), (0,0), (0,0), \cdots, (0,0) \}, \quad 
\{(1,1), (1,1), (0,1), (0,0), \cdots, (0,0) \}, \\
&\{(1,1), (0,1), (0,0), (0,0), \cdots, (0,0) \}. 
\end{align*}
}
By (b), the first one is $G_g$-equivalent to 
$\{(0,0), (0,0), (1,1), (0,0), \cdots, (0,0) \}$. 
Here, we observe actions of $G_g$ on the first and the second blocks, 
\begin{align*}
& \{(0,1), (1,1), \cdots \} \underset{C_1}{\longrightarrow}
\{(1,1), (1,1), \cdots \} \underset{\inv{DB_4} \cdot C_3 \cdot C_2}{\longrightarrow} 
\{(0,0), (0,1), \cdots \}, \\
& \{(0,1), (0,1), \cdots \} \underset{C_1}{\longrightarrow} 
\{(1,1), (0,1), \cdots \} \underset{C_3 C_2}{\longrightarrow} 
\{(0,0), (1,1), \cdots \}. 
\end{align*}
By the above observation, we see, 
\begin{align*}
\{(0,1), (1,1), (0,0), \cdots, (0,0) \} &\Ggeq 
\{(1,1), (1,1), (0,0), \cdots, (0,0) \}  \\
&\Ggeq\{(0,0), (0,1), (0,0), \cdots, (0,0) \}, \\
\{(0,1), (1,1), (0,1), \cdots, (0,0) \} &\Ggeq 
\{(1,1), (1,1), (0,1), \cdots, (0,0) \} \\
&\Ggeq \{(0,0), (0,1), (0,1), \cdots, (0,0) \} \\ 
&\Ggeq \{(0,0), (0,1), (0,0), \cdots, (0,0) \}, \\
\{(0,1), (0,1), (0,0), \cdots, (0,0) \} &\Ggeq 
\{(1,1), (0,1), (0,0), \cdots, (0,0) \} \\
&\Ggeq \{(0,0), (1,1), (0,0), \cdots, (0,0) \}. 
\end{align*}
Hence, we showed that any $x$ is $G_g$-equivalent to the elements 
listed in the statement of this Lemma. 
\end{proof}
Since 
\begin{align*}
&T^2_{\{ (0,1),(0,0), \cdots, (0,0) \}} = D_2, \quad
T^2_{\{ (1,1),(0,0), \cdots, (0,0) \}} = (C_1 C_2 C_1^{-1})^2, \\
&T^2_{\{ (0,0),(1,1),(0,0), \cdots, (0,0) \}} = (Y^*_2)^2, \quad
T^2_{\{ (0,0),(0,1),(0,0), \cdots, (0,0) \}} = D_4, \\
&T^2_{\{ (0,0),(0,0),(1,1), (0,0), \cdots, (0,0) \}} = (Y^*_4)^2, \quad
T^2_{\{ (0,0), \cdots, (0,0) \}} = id, 
\end{align*}
these are elements of $G_g$. 
By this fact and Lemma \ref{lem:Torelli}, 
Lemma \ref{lem:level2-congruence} is proved. 
\subsection{Step 3 for the case where $g \geq 3$}\label{subsec:Step3}
%
%
As in the previous subsection, let $\Phi_2 \co \mathcal{M}_g \to \Symp (2g, {\Bbb
Z}_2)$  be the natural homomorphism. 
Let $q_1 \co H_1(\Sigma_g; {\Bbb Z}_2) \to {\Bbb Z}_2$ be 
the quadratic form associated with the intersection form $(,)_2$ of 
$H_1(\Sigma_g; {\Bbb Z}_2)$ which satisfies, 
for the basis $x_i, y_i$ of $H_1(\Sigma_g; {\Bbb Z}_2)$ 
indicated on Figure \ref{fig:basis}, 
$q_1 (x_1) = q_1 (y_1) = 1$, and $q_1 (x_i) = q_1 (y_i) = 0$ when $i \not= 1$. 
We define $\Ortho_{q_1} (2g, {\Bbb Z}_2)$ $=$ 
$\{ \phi \in \Aut (H_1 (\Sigma_g; {\Bbb Z}_2)) | 
q_1 ( \phi(x)) = q_1 (x) \text{ for any } x \in 
H_1 (\Sigma_g; {\Bbb Z}_2) \}$, 
then $\oddspin{g}$ $=$ $\Phi_2^{-1}(\Ortho_{q_1}(2g, {\Bbb Z}_2))$. 
Because of Lemma \ref{lem:level2-congruence}, if we show 
$\Phi_2(G_g) = \Ortho_{q_1}(2g, {\Bbb Z}_2)$, then 
$G_g = \oddspin{g}$ follows.  

For any $z$ $\in$ $H_1 (\Sigma_g; {\Bbb Z}_2)$ such that 
$q_1(z)=1$, we define 
${\Bbb T}_z (x)= x + (z,x)_2 \; z$. 
Then ${\Bbb T}_z$ is an element of $\Ortho_{q_1} (2g, {\Bbb Z}_2)$, and 
we call this a {\it ${\Bbb Z}_2$-transvection about $z$\/}. 
Dieudonn\'e \cite{Dieudonne} showed the following Theorem 
(see also \cite[Chap.14]{Grove}). 
\begin{thm}\label{thm:Dieudonne}{\rm\cite[Proposition 14 on p.42]{Dieudonne}}\qua
%
When $g \geq 3$, $\Ortho_{q_1}(2g, {\Bbb Z}_2)$ is generated by 
${\Bbb Z}_2$-transvections. 
\end{thm}    
Let $\Lambda_g$ be the set of $z$ of $H_1 (\Sigma_g; {\Bbb Z}_2)$ 
such that $q (z) = 1$. 
For any elements $z_1$ and $z_2$ of $\Lambda_g$, we define 
$z_1 \square z_2 = z_1 + (z_2, z_1)_2 \; z_2$. 
Here, we remark that ${\Bbb T}_{z_1}^2 = \mathrm{id}$, 
${\Bbb T}_{z_2} {\Bbb T}_{z_1} {\Bbb T}_{z_2}^{-1}$ $=$ 
${\Bbb T}_{z_1 \square z_2}$ and $(z_1 \square z_2) \square z_2 = z_1$. 
An element 
$\epsilon_1 x_1 + \delta_1 y_1 + \cdots + \epsilon_g x_g + \delta_g y_g$ 
of $H_1 (\Sigma_g; {\Bbb Z}_2)$ is denoted by 
$[(\epsilon_1, \delta_1), \cdots, (\epsilon_g, \delta_g)]$, 
and each $(\epsilon_i, \delta_i)$ is called the $i$-th block. 
We remark that $q([(\epsilon_1, \delta_1), \cdots, (\epsilon_g, \delta_g)]) 
= (\epsilon_1 + \delta_1 + \epsilon_1 \delta_1) + \epsilon_2 \delta_2 + 
\cdots \epsilon_g \delta_g$. 
\begin{lem}\label{lem:generator-lambda-g}
%
Under the operation $\square$, $\Lambda_g$ is generated by 
$x_1$, $y_1$, $x_1+x_2$, 
$x_i + y_i$ $(2 \leq i \leq g)$, $x_i + y_i + x_{i+1}$ $(2 \leq i \leq g-1)$, 
and $x_i+x_{i+1}+y_{i+1}$ $(2 \leq i \leq g-1)$. 
\end{lem}
\begin{proof}
For an element $[(\epsilon_1, \delta_1), \cdots, (\epsilon_g, \delta_g)]$ of 
$H_1 (\Sigma_g; {\Bbb Z}_2)$, let the $j$-th block be the right most block 
which is $(1,1)$. 
When $j \geq 3$, there exist 4 cases of the combination of the $(j-1)$-st
block  and the $j$-th block: 
$[\cdots, (1,1), (1,1), \cdots]$, $[\cdots, (0,0), (1,1), \cdots]$, 
$[\cdots, (0,1), (1,1), \cdots]$, $[\cdots, (1,0), (1,1), \cdots]$. 
In each case, we can reduce $j$ at least 1. In fact, 
%
\begin{align*}
[\cdots, (1,1), (1,1), \cdots] &\square (x_{j-1} + x_j + y_j) 
=[\cdots, (0,1),(0,0),\cdots], \\
[\cdots, (0,0), (1,1), \cdots] &\square (x_{j-1} + y_{j-1} + x_j) 
=[\cdots, (1,1),(0,1),\cdots], \\
[\cdots, (0,1), (1,1), \cdots] &\square (x_{j-1} + x_j + y_j) 
=[\cdots, (1,1),(0,0),\cdots], \\
([\cdots, (1,0), (1,1), \cdots] &\square (x_{j-1} + y_{j-1})) 
\square (x_{j-1} + x_j + y_j) \\
&=[\cdots, (1,1),(0,0),\cdots]. 
\end{align*}
%
When $j=2$, 
since $q([(\epsilon_1, \delta_1), \cdots, (\epsilon_g, \delta_g)]) =1$, 
$[(\epsilon_1, \delta_1), \cdots, (\epsilon_g, \delta_g)]$ must be 
$[(0,0),(1,1),\cdots]$. 
Because of an equation 
\begin{equation*}
([(0,0),(1,1),\cdots] \square (x_1+x_2)) \square y_1 = [(1,1),(0,1),\cdots], 
\end{equation*}
we can reduce $j$ to $1$. 
When $j=1$, if every $i$-th ($i \geq 2$) block is $(0,0)$, then 
it is $x_1 + y_1$, which is equal to $x_1 \square y_1$. 
If there exist at least one of 
the $i$-th ($i \geq 2$) blocks which are $(1,0)$ or $(0,1)$, then,  
%
\begin{align*}
&[\cdots, (0,0), \overset{i}{(1,0)}, \cdots] \square (x_{i-1} + x_i + y_i) 
=[\cdots, (1,0), (0,1), \cdots], \\
&[\cdots, (1,0), \overset{i}{(0,0)}, \cdots] \square (x_{i-1} + y_{i-1} +x_i) 
=[\cdots, (0,1), (1,0), \cdots], \\
&[\cdots, (0,0), \overset{i}{(0,1)}, \cdots] \square (x_{i-1} + x_i + y_i) 
=[\cdots, (1,0), (1,0), \cdots], \\
&[\cdots, (0,1), \overset{i}{(0,0)}, \cdots] \square (x_{i-1} + y_{i-1}+x_i ) 
=[\cdots, (1,0), (1,0), \cdots].
\end{align*}
%
Therefore, we can alter this to an element, 
each $i$-th ($i \geq 2$) block of which is $(1,0)$ or $(0,1)$. 
If the $i$-th block of this is $(0,1)$, then 
$$
[\cdots, (0,1), \cdots] \square (x_i + y_i) = [\cdots, (1,0), \cdots].
$$ 
Therefore, it suffices to consider the case where 
the first block is $(1,1)$ and other blocks are $(1,0)$. 
In this case, 
$$
([\cdots, (1,0), (1,0)] \square (x_{g-1}+y_{g-1}+x_g) )
\square (x_{g-1} + y_{g-1}) = [\cdots, (1,0), (0,0)]. 
$$
By applying the same operation repeatedly, we get 
$[(1,1),(1,0),(0,0), \cdots, $ $(0,0)]$, 
which is equal to $y_1 \square (x_1 + x_2)$. 
\end{proof}
This lemma and Theorem \ref{thm:Dieudonne} shows that  
\begin{cor}\label{cor:generator-o-2g}
%
$\Ortho_{q_1} (2g, {\Bbb Z}_2)$ is generated by 
${\Bbb T}_{x_1}$, ${\Bbb T}_{y_1}$, ${\Bbb T}_{x_1 + x_2}$, 
${\Bbb T}_{x_i + y_i}$ $(2 \leq i \leq g)$, 
${\Bbb T}_{x_i + y_i +x_{i+1}}$ $(2 \leq i \leq g-1)$, and 
${\Bbb T}_{x_i + x_{i+1} + y_{i+1}}$ $(2 \leq i \leq g-1)$.
\qed
\end{cor}
Since $G_g$ is a subgroup of $\oddspin{g}$,  
$\Phi_2 (G_g) \subset \Ortho_{q_1} (2g, {\Bbb Z}_2)$. 
On the other hand, the fact that 
$\Phi_2(C_1) = \mathbb{T}_{x_1}$, $\Phi_2(C_2) = \mathbb{T}_{y_1}$, 
$\Phi_2(C_3) = \mathbb{T}_{x_1 + x_2}$, 
$\Phi_2(X_{2i}) = {\Bbb T}_{x_i + y_i + x_{i+1}}$ ($2 \leq i \leq g-1$), 
$\Phi_2(X_{2i+1}) = {\Bbb T}_{x_i + x_{i+1} + y_{i+1}}$ ($2 \leq i \leq g-1$), 
$\Phi_2(Y_{2j}) = {\Bbb T}_{x_j + y_j}$ ($2 \leq j \leq g-1$), 
$\Phi_2(X_{2g}) = {\Bbb T}_{x_g + y_g}$, 
and Corollary \ref{cor:generator-o-2g}, 
show $\Phi_2 (G_g) \supset \Ortho_{q_1} (2g, {\Bbb Z}_2)$. 
Therefore we proved that $\oddspin{g} = G_g$ when $g \geq 3$. 
\subsection{Genus 2 case: Reidemeister-Schreier method}\label{subsec:genus2}
%
%
Birman and Hilden showed the following Theorem. 
\begin{thm}\label{thm:genus2}{\rm\cite{Birman-Hilden}}\qua
%
$\mathcal{M}_2$ is generated by $C_1, C_2, C_3, C_4, C_5$ and its defining 
relations are: 
\par\noindent
(1) $C_i C_j = C_j C_i$, if $|i-j| \geq 2$, $i,j=1,2,3,4,5$, 
\par\noindent
(2) $C_i C_{i+1} C_i = C_{i+1} C_i C_{i+1}$, $i=1,2,3,4$, 
\par\noindent
(3) $(C_1 C_2 C_3 C_4 C_5)^6 = 1$, 
\par\noindent
(4) $(C_1 C_2 C_3 C_4 C_5 C_5 C_4 C_3 C_2 C_1)^2 = 1$, 
\par\noindent
(5) $C_1 C_2 C_3 C_4 C_5 C_5 C_4 C_3 C_2 C_1 \rightleftarrows C_i$, 
$i=1,2,3,4,5$, 
\par\noindent
where $\rightleftarrows$ means "commute with". 
\end{thm}
We call (1) (2) of the above relations {\it braid relations\/}. 
We will use the well-known method, called 
{\it the Reidemeister--Schreier method\/} \cite[\S 2.3]{MKS}, 
to show $\oddspin{2} \subset G_2$. 
We review (a part of) this method. 
\par
Let $G$ be a group generated by finite elements $g_1, \ldots, g_m$ and 
$H$ be a finite index subgroup of $G$. 
For two elements $a$, $b$ of $G$, we write $a \equiv b$ mod $H$ if 
there is an element $h$ of $H$ such that $a = hb$. 
A finite subset $S$ of $G$ is called {\it a coset representative system\/} 
for $G$ mod $H$, if, for each elements $g$ of $G$, there is only one element 
$\repre{g} \in S$ such that $g \equiv \repre{g}$ mod $H$. 
The set $\{ s g_i \repre{s g_i}^{-1} \ | \ i=1,\ldots,m,\ s \in S \}$ 
generates $H$. 
\par
\begin{figure}[ht!]\small
\centering
\SetLabels
(0*0.7) $[0,1,1,1]$ \\ (0.2*0.7) $[0,0,1,1]$ \\ (0.4*0.7) $[1,0,1,1]$ \\
(0.6*0.7) $[1,1,1,0]$ \\ (0.8*0.7) $[1,1,0,0]$ \\ (1*0.7) $[1,1,0,1]$ \\
(0.1*0.15) $C_1$ \\ (0.3*0.15) $C_2$ \\ (0.5*0.15) $C_3$ \\ 
(0.69*0.15) $C_4$ \\ (0.89*0.15) $C_5$ \\
\endSetLabels
\strut\AffixLabels{\includegraphics[height=1.5cm]{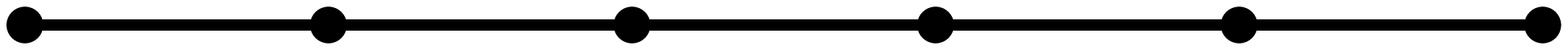}}
\nocolon
\caption{}
\label{fig:oddaction}
\end{figure}
For the sake of giving a coset representative system for $\mathcal{M}_2$ modulo 
$\oddspin{2}$, we will draw a graph $\Gamma$ which represents the action of 
$\mathcal{M}_2$ on the quadratic forms of $H_1(\Sigma_2; {\Bbb Z}_2)$ with 
Arf invariants $1$. 
Let $[\epsilon_1, \epsilon_2, \epsilon_3, \epsilon_4]$ denote the quadratic form 
$q'$ of $H_1(\Sigma_2; {\Bbb Z}_2)$ such that $q'(x_1)=\epsilon_1$, 
$q'(y_1)=\epsilon_2$, $q'(x_2)=\epsilon_3$, $q'(y_2)=\epsilon_4$. 
Each vertex of $\Gamma$ corresponds to a quadratic form. 
For each generator $C_i$ of $\mathcal{M}_2$, we denote its action on 
$H_1(\Sigma_2; {\Bbb Z}_2)$ by $(C_i)_*$. 
For the quadratic form $q'$ indicated by the symbol 
$[\epsilon_1, \epsilon_2, \epsilon_3, \epsilon_4]$, 
let $\delta_1 = q'((C_i)_* x_1)$,  $\delta_2 = q'((C_i)_* y_1)$, 
 $\delta_3 = q'((C_i)_* x_2)$,  and $\delta_4 = q'((C_i)_* y_2)$. 
Then, we connect two vertices, corresponding to 
$[\epsilon_1, \epsilon_2, \epsilon_3, \epsilon_4]$, 
$[\delta_1, \delta_2, \delta_3, \delta_4]$ respectively, 
by the edge with the letter $C_i$. 
We remark that this action is a right action. 
For simplicity, we omit the edge whose ends are the same vertex. 
As a result, we get a graph $\Gamma$ as in Figure \ref{fig:oddaction}. 
The words $S=\{1,\ C_5,\ C_4,\ C_4 C_3,\ C_4 C_3 C_2,\ C_4 C_3 C_2 C_1\}$, 
which correspond to the edge paths beginning 
from $[1,1,0,0]$ on $\Gamma$, define a coset representative system for 
$\mathcal{M}_2$ modulo $\oddspin{2}$. 
For each element $g$ of $\mathcal{M}_2$, we can give a $\repre{g}$ $\in S$ with 
using this graph. 
For example, say $g=C_2 C_4 C_5 C_3$, we follow an edge path assigned 
to this word which begins from $[1,1,0,0]$, 
(note that we read words from left to right) then we arrive at the vertex 
$[1,0,1,1]$. 
The element in $T$ which begins from $[1,1,0,0]$ and ends at $[1,0,1,1]$ 
is $C_4 C_3$. Hence, $\repre{C_2 C_4 C_5 C_3} = C_4 C_3$. 
We list in Table \ref{tab:generators} the set of generators 
$\{ s C_i \repre{s C_i}^{-1} \ | \ i=1,\ldots,5,\ s \in S\}$ of 
$\oddspin{g}$. 
In Table \ref{tab:generators}, 
vertical direction is a coset representative system 
$S$, horizontal direction is a set of generators 
$\{C_1,\ C_2,\ C_3,\ C_4,\ C_5\}$. 
\begin{table}
\caption{Generators of $\oddspin{2}$}
\label{tab:generators}
\begin{center}
\begin{tabular}{c|c c c}
 & $C_1$ & $C_2$ & $C_3$\\
\hline
$1$ & $1$ & $1$ & $1$ \\
$C_5$ & $C_1$ & $C_2$ & $C_3$ \\
$C_4$ & $C_1$ & $C_2$ & $1$\\
$C_4 C_3$ & $C_1$ & $1$ & $C_3^{-1} D_4 C_3$ \\ 
$C_4 C_3 C_2$ & $1$ & $C_2^{-1} C_3^{-1} D_4 C_3 C_2$ & $C_2$ \\
$C_4 C_3 C_2 C_1$ & $C_1^{-1} C_2^{-1} C_3^{-1} D_4 C_3 C_2 C_1$ & $C_1$ & $C_2$\\ 
\end{tabular}
\begin{tabular}{c|c c}
&$C_4$ & $C_5$ \\
\hline
$1$ & $1$ & $1$ \\
$C_5$ & $1$ & $D_5$ \\
$C_4$ & $D_4$ & $D_5^{-1} X_4 D_5$ \\
$C_4 C_3$ & $C_3$ & $D_5^{-1} X_4 D_5$ \\ 
$C_4 C_3 C_2$ &$C_3$ & $D_5^{-1} X_4 D_5$ \\
$C_4 C_3 C_2 C_1$ & $C_3$ & $D_5^{-1} X_4 D_5$ \\ 
\end{tabular}
\end{center}
\end{table} 
We can check this table by Figure \ref{fig:oddaction} and braid relations. 
For example, 
\begin{align*}
&C_4 C_3 C_2 C_1 \cdot C_2 \repre{ C_4 C_3 C_2 C_1 \cdot C_2 }^{-1} = 
C_4 C_3 C_2 C_1 C_2 (C_4 C_3 C_2 C_1)^{-1} \\
&= C_4 C_3 C_2 C_1 C_2 C_1^{-1} C_2^{-1} C_3^{-1} C_4^{-1} = 
C_4 C_3 C_2 C_2^{-1} C_1 C_2 C_2^{-1} C_3^{-1} C_4^{-1} \\
&= C_4 C_3 C_1 C_3^{-1} C_4^{-1} = C_1. 
\end{align*}
This table shows that 
$\oddspin{2} \subset G_2$ . 

\section{Proof of Theorem \ref{thm:main}}\label{sec:extension}

\begin{figure}[ht!]\small
\centering
\SetLabels
(0.14*0.9) $c_2$ \\ (0.17*0.5) $c_1$ \\ (0.37*0.7) $c_3$ \\ 
(0.50*0.55) $c_4$ \\ (0.46*0.93) $b_4$ \\ (0.67*0.53) $c_{2g-2}$ \\
(0.76*0.67) $c_{2g-1}$ \\ (0.67*0.94) $b_{2g-2}$ \\ 
(0.86*0.53) $c_{2g}$ \\ \L(0.94*0.67) $c_{2g+1}$ \\
\endSetLabels
\strut\AffixLabels{\includegraphics[height=4.5cm]{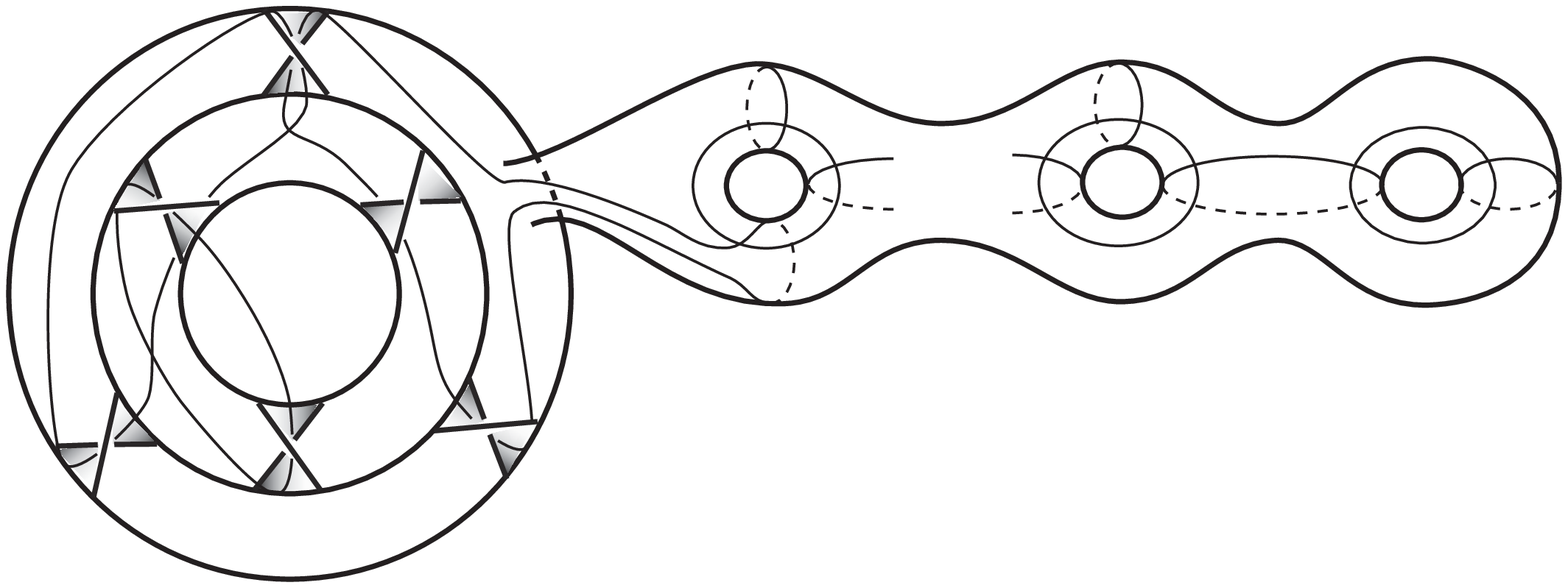}}
\nocolon
\caption{}
\label{fig:connected-sum}
\end{figure}
We embed $H_{g-1}$ standardly in $S^3 = \partial D_4$ such that there is 
a 2-sphere separating $F_{3,3}$ and $H_{g-1}$, and make a connected sum 
$F_{3,3} \# \partial H_{g-1}$ as indicated in 
Figure \ref{fig:connected-sum}. 
Then, we can see 
$(\ComplexPlane, K_3 \# \Sigma_{g-1})$ $=$ 
$(\ComplexPlane, (F_{3,3} \# \partial H_{g-1}) \cup D_3)$, 
where $K_3$ is the non-singular plane curve of degree $3$ and 
$D_3$ is parallel three disks which is used to construct $K_3$ 
in \S \ref{sec:curve}. 
We identify $K_3 \# \Sigma_{g-1}$ with $\Sigma_g$ so that simple closed curves with 
the same symbol are identified. 
Then $q_{K_3 \# \Sigma_{g-1}} = q_1$. 
We will show that each elements of $\spin{g}{q_{K_3 \# \Sigma_{g-1}}}$ $=$ 
$\spin{g}{q_1}$ is extendable. 

Each regular neighborhood of $c_1$, $c_2$, $c_3$, $C_{i+1} (c_i)$ 
($4 \leq i \leq 2g$), and $C_{2j} (b_{2j})$ ($2 \leq j \leq g-1$) is Hopf band. 
Therefore, by Proposition \ref{prop:Hopf-twist}, 
$C_1$, $C_2$, $C_3$, $C_{i+1} C_i \inv{C_{i+1}}$ ($4 \leq i \leq 2g$), and 
$C_{2j} B_{2j} \inv{C_{2j}}$ ($2 \leq j \leq g-1$) are elements of 
$\mathcal{E} (\ComplexPlane, K_3 \# \Sigma_{g-1})$. 
Each regular neighborhood of $c_i$ ($4 \leq i \leq 2g+1$), $b_{2j}$ 
($2 \leq i \leq g-1$) is an annulus standardly embedded in $S^3 = \partial D^4$. 
We can deform this annulus as indicated in Figure \ref{fig:bandmove}. 
Therefore, $C_i^2$ ($4 \leq i \leq 2g+1$), $B_{2j}^2$ ($2 \leq j \leq g-1$) are 
elements of $\mathcal{E}(\ComplexPlane, K_3 \# \Sigma_{g-1})$. 
Finally, the extendability of $B_4 C_5 C_7 \ldots C_{2g+1}$ follows from the proof 
of Lemma 2.2 in \cite{Hirose2}. 
Therefore, we showed $\spin{g}{q_{K_3 \# \Sigma_{g-1}}}$ $\subset$ 
$\mathcal{E}(\ComplexPlane, K_3 \# \Sigma_{g-1})$. 
On the other hand, by the definition of the Rokhlin quadratic form 
$q_{K_3 \# \Sigma_{g-1}}$, we see $\mathcal{E}(\ComplexPlane, K_3 \# \Sigma_{g-1})$
$\subset$ $\spin{g}{q_{K_3 \# \Sigma_{g-1}}}$. 
Theorem \ref{thm:main} follows. 

\subsection*{Acknowledgments}
The author would like to express his gratitude to Professors 
Masaharu Ishikawa, Masahico Saito, and Akira Yasuhara 
for fruitful discussions and comments. 
The author would also like to thank the referee, 
whose comments and corrections improved the paper.  
This research was partially supported by Grant-in-Aid for 
Encouragement of Young Scientists (No. 16740038), 
Ministry of Education, Culture, Sports, Science and Technology, Japan.

%
\Addresses \recd
\end{document}